\documentclass[letterpaper,12pt,]{amsart}
\usepackage{preambleFL}

\usepackage{enumerate}
\usepackage[paper=letterpaper,right=3cm,left=3cm,top=2.5cm,bottom=2.5cm]{geometry}
\usepackage{stmaryrd} 
\usepackage{makecell}

 \newtheorem*{thm*}{Theorem}
 \newtheorem*{defn*}{Definition}

\newcommand{\universal}{\zeta_n}

\newcommand{\moore}{\ensuremath{\mathbf{\Lambda}}}

\newcommand{\topG}{\ensuremath{\mathrm{Top}_G}}
\newcommand{\ele}{\ensuremath{\nu}}
\newcommand{\elm}{\ensuremath{\mu}}
\newcommand{\subgroup}{\ensuremath{\subset}}

\newif\ifcomment
\commentfalse
\newif\ifall
\allfalse

\begin{document}

\title[Equivariant bundles]{Notes on equivariant bundles}
\author{Foling Zou}

 \address{Department of Mathematics, University of Michigan, Ann Arbor, MI 48109 United States}
 \email{folingz@umich.edu}

\maketitle
\begin{abstract}
We compare two notions of $G$-fiber bundles and $G$-principal bundles in the
literature, with an aim to clarify early results in equivariant bundle
theory that are needed in current work of equivariant algebraic topology.
We also give proofs of some equivariant generalizations of well-known non-equivariant
results involving the classifying space.
\end{abstract}

\tableofcontents

\section{Introduction}
Non-equivariantly, fiber bundles and principal bundles are closely
related. Namely, fixing
a compact Lie group $\Pi$ and a space $F$ with an effective $\Pi$-action,
one can make sense of a fiber bundle with fiber $F$ to have structure group $\Pi$,
and there is a structure theorem providing an equivalence of categories between such fiber bundles
and principal $\Pi$-bundles. One key idea involved is the data of admissible maps of a
fiber bundle $p:E \to B$ with fiber $F$, which are specified homeomorphisms
$\psi: F \cong p^{-1}(b)$ for $b \in B$ that come from the local trivializations.

Equivariantly, let $G$ and $\Pi$ be compact Lie groups, with $G$ being the ambient action group and $\Pi$ being the
structure group. It is common to assume $G$ to be a compact Lie
  group in $G$-equivariant homotopy theory for several reasons. Firstly, the orbit category
  of a compact Lie group $G$ is more controlled as it has a
  discrete skeleton. More precisely, the orbit category has objects $G/H$ for closed subgroups $H
  \subset G$, and two closed subgroups of $G$ being sufficiently close will
  imply that one is conjugate to a subgroup of the other by
  \cite{MontgometryZippin}. If one works with a general topological group $G$, the isotropy subgroups of a $G$-space may vary continuously.
  Secondly, nice $G$-spaces such as
  smooth $G$-manifolds are known to allow $G$-equivariant triangulation
  by \cite{IllmanCLie}. One can work with more general topological
  groups $\Pi$ with careful point set considerations, but we refrain from doing so in this
  paper.

To obtain a structure theorem relating fiber bundles and principal bundles equivariantly, we need to answer the following two
questions: What does it mean for a $G$-fiber bundle with fiber $F$ to have
structure group $\Pi$? What is an equivariant principal $\Pi$-bundle?

\autoref{sec:defin-equiv-bundl} is devoted to answering these two questions and
establishing the structure theorems.
The forthright guess already works for some examples including equivariant
vector bundles. We assume that the $\Pi$-action on $F$ is
effective throughout.

\begin{defn}(\autoref{defn:Gvector})
  \label{defn:intro1}
  Let $F$ be a space with $\Pi$-action. A $G$-fiber bundle with fiber $F$ and structure
  group $\Pi$ is a map $p:E \to B$ such that the following statements hold:
  \begin{enumerate}
    \item The map $p$ is a non-equivariant fiber bundle with fiber $F$ and
      structure group~$\Pi$;
    \item Both $E$ and $B$ are $G$-spaces and $p$ is $G$-equivariant;
    \item \label{item:Gvector3} The $G$-action is given by morphisms of bundles with structure group $\Pi$.
    \end{enumerate}
  \end{defn}

  Tom~Dieck \cite{TD} generalized this definition to a twisted version and
   Lewis--May--Steinberger \cite[IV1]{LMS86} introduced the following further generalization. We fix an
   extension of compact Lie groups $1 \to \Pi \to \Gamma \to G \to 1$ as
   data.
   
   \begin{defn}(\autoref{defn:Gfiber})
   \label{defn:intro2}
   Let $F$ be a space with $\Gamma$-action.
   A $G$-fiber bundle with fiber $F$, structure group $\Pi$ and total group $\Gamma$
   is a map $p:E \to B$ such that the following statements hold:
  \begin{enumerate}
    \item The map $p$ is a non-equivariant fiber bundle with fiber $F$ and structure group~$\Pi$;
    \item Both $E$ and $B$ are $G$-spaces and $p$ is a $G$-map;
    \item \label{item:Gfiber-admissible}
      For any  $g \in G$ and admissible maps
      $\psi: F \to F_b$ and $\zeta: F \to F_{gb}$, the composite      
      \begin{equation*}
        F \overset{\psi}{ \to } F_b \overset{g}{ \to } F_{gb} \overset{\zeta^{-1}}{ \to } F
      \end{equation*}
      is a lift $y \in \Gamma$ of $g \in G$.
    \end{enumerate}
  \end{defn}
  In terms of admissible maps, \autoref{defn:intro1} requires that for any admissible
map $\psi: F \to E$ and $g \in G$, the composite $g \psi$ is also admissible;
while \autoref{defn:intro2} requires that there is a lift $y \in \Gamma$ of $g$
such that $g \psi y^{-1}$ is admissible.

The following gives an interesting example of the generalization
.
  \begin{exmp}(Examples~\ref{eg:Real}, \ref{example:real} and \ref{exmp:RealNot})
  Atiyah \cite{Atiyah} introduced the notion of Real vector bundles: it is a
  complex vector bundle with $C_2$-action such that the non-trivial element of
  $C_2$ acts anti-complex-linearly.
  It is a $C_2$-fiber bundle with fiber $\bC^n$, structure
  group $U(n)$ and total group $\Gamma = U(n) \rtimes_{\alpha} C_2$ in the sense
  of \autoref{defn:intro2}, where $\alpha: C_2 \to
  \mathrm{Aut}(U(n))$ sends the non-trivial element of $C_2$ to the entry-wise complex-conjugation
  of $U(n)$. But it is \emph{not} a $C_2$-fiber bundle with fiber $\bC^n$, structure
  group $U(n)$ in the sense of \autoref{defn:intro1}.
\end{exmp}

Note that in \autoref{defn:intro1}, the fiber $F$ only has a $\Pi$-action, but
in \autoref{defn:intro2}, the fiber $F$ has an action of the total group $\Gamma$. 
\autoref{defn:intro1} becomes a special case of \autoref{defn:intro2} by setting $\Gamma = \Pi
\times G$ and imposing the trivial
$G$-action  on $F$ (see \autoref{prop:bundle-action1}).
Tom~Dieck's definition, which is not given above, corresponds to the case of split extensions $\Gamma = \Pi
   \rtimes_{\alpha} G$ in \autoref{defn:intro2}. It allows the
   $G$-action on the bundle to have a preassigned twisting in the structure
   group, such as in the Real vector bundles. Conceptually, the most general
   \autoref{defn:intro2} sees a more general twisting specified by the
   group extension  $1 \to \Pi \to \Gamma \to G \to 1$. 

The following are companion definitions of equivariant principal bundles.

\begin{defn}(\autoref{defn:Gprincipal})
  A principal $G$-$\Pi$-bundle is a map ${p:P \to B}$ such that the following statements hold:
    \begin{enumerate}
      \item The map $p$ is a non-equivariant principal $\Pi$-bundle;
      \item Both $P$ and $B$ are $G$-spaces and $p$ is $G$-equivariant;
      \item The actions of $G$ and $\Pi$ commute on $P$.
    \end{enumerate}
\end{defn}

\begin{defn}(\autoref{defn:GGammaprincipal})
  Let $1 \to \Pi \to \Gamma \to G \to 1$ be an extension of compact Lie groups. 
  A principal $(\Pi;\Gamma)$-bundle is a map $p:P \to B$ such that the
  following statements hold:
    \begin{enumerate}
      \item The map $p$ is a non-equivariant principal $\Pi$-bundle;
      \item The space $P$ is a $\Gamma$-space; $B$ is a $G$-space. Viewing $B$ as a
 $\Gamma$-space by pulling back the action, the map $p$ is $\Gamma$-equivariant.
    \end{enumerate}
  \end{defn}

  The first definition is a special case of the second definition by setting $\Gamma
  = \Pi \times G$.

  \medskip
  
There are structure theorems relating equivariant fiber bundles and equivariant
principal bundles (Theorems~\ref{thm:G-structure-1} and  \ref{thm:G-structure-2}).
For any $G$-space $B$, there is an equivalence of categories between

\smallskip
    \begin{tabular}{ccc}
      \makecell{$G$-fiber bundles over $B$ with fiber $F$
      \\ and structure group $\Pi$}
      & $\leftrightarrow$
      & principal $G$-$\Pi$-bundles over $B$. \\
      \makecell{$G$-fiber bundles over $B$ with  fiber $F$
      \\ structure group $\Pi$ and total group $\Gamma$}
      & $\leftrightarrow$
      & principal $(\Pi;\Gamma)$-bundles over $B$. 
    \end{tabular}
    \smallskip
    
    \noindent  In the first equivalence,
    $F$ is any $\Pi$-effective space; in the second
 equivalence, $F$ is any $\Pi$-effective $\Gamma$-space.
 
  As an example, we study the $V$-framing bundle $\mathrm{Fr}_V(E)$ of a
  $G$-$n$-vector bundle $E$, where $V$ is a
  $G$-representation given by $\beta: G \to O(n)$. It turns out 
  $\mathrm{Fr}_{\bR^n}(E) \cong \mathrm{Fr}_V(E)$ as $(O(n); O(n) \times G) \cong (O(n);O(V) \rtimes
  G)$-principal bundles in the sense of \autoref{defn:intro2}, hinting that the
  $V$-framing bundle may not be an interesting notion. However, they can be
  given different canonical $G$-actions and
  the $G$-action on $\mathrm{Fr}_V(E)$ can be
  identified with the
  $\Lambda_{\beta}$-action on $\mathrm{Fr}_{\bR^n}(E)$ for some subgroup
  $\Lambda_{\beta} \subset O(n) \times G$. This is in
  \autoref{subsec:split-ext}.

  \medskip
There exists a universal principal $(\Pi; \Gamma)$-bundle $E(\Pi; \Gamma) \to B(\Pi; \Gamma)$. It is universal in the sense that
 there is a bijection of sets between
\{equivalence classes of principal $(\Pi;\Gamma)$-bundles over $B$\} and  \{$G$-homotopy classes of $G$-maps $B \to B(\Pi;\Gamma)$\} for
any paracompact $G$-space $B$. Thus, $B(\Pi; \Gamma)$ is called the
classifying space of principal $(\Pi; \Gamma)$-bundles.
In the case $\Gamma = \Pi \times G$, we also denote the universal principal $G$-$\Pi$-bundle by $E_G\Pi \to  B_G\Pi.$

The universal principal bundle can be constructed using homotopy theory techniques.
Non-equivariantly, one can construct $E \Pi$, a contractible space with free $\Pi$-action, and the universal principal $\Pi$-bundle is modeled
by $E\Pi \to E \Pi/\Pi$. Equivariantly, a family  $\mathscr{F}$ of subgroups of $\Gamma$ is
a collection of subgroups  that is non-empty and closed under
subgroups and conjugations. For each family, there exists a $\Gamma$-space $E \mathscr{F}$
with the property that $(E \mathscr{F})^{\Lambda} \simeq \begin{cases}
  \varnothing & \Lambda \subset \Gamma \text{ and } \Lambda \not\in \mathscr{F} \\
  * & \Lambda \subset \Gamma \text{ and } \Lambda \in \mathscr{F}
\end{cases}$. We take the family to be $\mathscr{F}  = \mathscr{F}(\Pi) :=
\{\Lambda \subset \Gamma | \Lambda \cap \Pi = e \}$. Then the universal
principal $(\Pi;\Gamma)$-bundle is modeled by $E\mathscr{F}(\Pi) \to E
\mathscr{F}(\Pi) / \Pi$ (\autoref{thm:universalbundle}). Note that to recover the
non-equivariant case, we can take $G = e$, so that $\Gamma = \Pi$ and $\mathscr{F}(\Pi) = e$, and
the fixed-point properties of $E \mathscr{F}(\Pi)$ coincide with the
defining properties
of $E \Pi$. More details are in \autoref{sec:univ-bundle}.

 \begin{rem}
\label{rem:univ-topo}
L\"uck--Uribe worked with general topological groups $G$ and $\Pi$ and
those principal $G$-$\Pi$-bundles
  such that the isotropy subgroups of the total space are in $\mathcal{R}$ for a
  prescribed family of subgroups of $\Gamma = \Pi \times G$.
  This family needs to satisfy conditions in \cite[Definition 6.1]{LU14}.
  The universal bundle of such principal bundles can be modeled by $E \mathcal{R} \to E \mathcal{R}/\Pi$
  (\cite[Theorem 11.5]{LU14}). 
  In our case, $\mathcal{R} =  \mathscr{F}(\Pi)$,
  and we will not make use of a general $\mathcal{R}$.
  To translate notations, their $\Gamma$ is our $G$ and their $G$ is our $\Pi$.
 \end{rem}

  \medskip
   In \autoref{sec:fixed-point-theorems}, we study the fixed points of a principal $G$-$\Pi$-bundle  $p: P
   \to B$. Let $H \subgroup G$ be a subgroup and use $\mathrm{Rep}(H,\Pi)$ to
   denote group homomorphisms $\rho: H \to \Pi$ up to $\Pi$-conjugation.
Each component $B_0$ of $B^H$ has an associated homomorphism $[\rho] \in
\mathrm{Rep}(H, \Pi)$. The $[\rho]$ is determined by the fixed-point behavior of the total space:
Let $\Lambda_{\rho} \subgroup \Pi \times G$ denote the graph of $\rho$,
then $\{\rho: H \to \Pi | \big(p^{-1}(B_0)\big)^{\Lambda_{\rho}} \not=
\varnothing \}$ holds for exactly one $\Pi$-conjugate class of homomorphisms.
Furthermore, the non-equivariant principal $\Pi$-bundle $p^{-1}(B_0) \to B_0$  has a reduction
of the structure group from $\Pi$ to a subgroup $Z_{\Pi}(\rho) \subgroup \Pi$ (\autoref{thm:LM}).
We apply this theorem to obtain a comparison  of principal $G$-$\Pi$-bundles. In a map between principal
$G$-$\Pi$-bundles, if $\bar{f}$, the map of the total spaces, is an
$\mathscr{F}(\Pi)$-equivalence, then $f$, the map of the base spaces, is a $G$-equivalence
(\autoref{lem:comparefib}).

\medskip
In \autoref{sec:fiber-univ-bundle}, we study the loop space $\Omega_bB_G\Pi$ of $B_G\Pi$ based at a $G$-fixed point $b$.
As $(B_G\Pi)^G$ is not connected in general, the $G$-homotopy type of $\Omega_bB_G\Pi$ depends on the
choice of $b$. Our greatest interest is in the case $\Pi=O(n)$, and it works the same for general
$\Pi$ as discussed in \autoref{rem:general-Pi}.
Note that $B_GO(n)$ classifies $G$-$n$-vector bundles and
a homomorphism $\rho: G \to O(n)$ gives an $n$-dimensional $G$-representation $V$.
Suppose $b \in B_GO(n)$ is in the component indexed by $[V]$. In \autoref{thm:BGO(n)2},
we show that there is a $G$-homotopy equivalence $\Omega_b B_GO(n) \simeq O(V)$, where $O(V)$ is the
isometric self maps of $V$ with $G$-action by conjugation. Later in
\autoref{cor:monoidmap}, we upgrade the $G$-equivalence to one compatible with
the monoid structure: $\Omega_b B_GO(n)$ has the structure of a
$G$-$A_{\infty}$-monoid via concatenation of loops. There is a zig-zag of
equivalence of $G$-$A_{\infty}$-monoids $\Omega_bB_GO(n) \simeq O(V)$.

In \autoref{sec:gauge-group}, we study the space of bundle maps to
  the equivariant universal bundle. It is an
equivariant principal bundle with a non-trivial extension in the total group.
Let $p: P \to B$ be a principal $G$-$O(n)$-bundle, $\Pi= \mathrm{Aut}_B(P)$
be the topological group of automorphisms
of $P$ over $B$ and $\mathrm{Hom}(P, E_GO(n))$ be the space of non-equivariant principal
$O(n)$-bundle maps. We have that $G$ acts on $\mathrm{Hom}(P, E_GO(n))$ and $\Pi$ by
conjugation, and that $\Pi$ acts on  $\mathrm{Hom}(P, E_GO(n))$ by precomposition.
This gives a $(\Gamma =\Pi \rtimes G)$-action on $\mathrm{Hom}(P, E_GO(n))$.
In \autoref{thm:MapToUniversal}, we show that
\begin{equation*}
\pi: \mathrm{Hom}(P, E_GO(n)) \to \mathrm{Map}_p(B,B_GO(n))
\end{equation*}
is a principal  $(\mathrm{Aut}_{B}P; \mathrm{Aut}_BP \rtimes G)$-bundle and $\mathrm{Hom}(P, E_GO(n)) \simeq E \mathscr{F}$ for the family $\mathscr{F} = \{\Lambda \subgroup  \Gamma \text{ such
  that } \Lambda \cap \mathrm{Aut}_BP = e\}$.
Here, $\pi$ sends a bundle map to its map of base spaces and $\mathrm{Map}_p(B,B_GO(n))$ is the
image of $\pi$ in $\mathrm{Map}(B,B_GO(n))$.
We conjecture that $\pi$ is the universal principal  $(\mathrm{Aut}_{B}P;
\mathrm{Aut}_BP \rtimes G)$-bundle. The issue is that $\mathrm{Aut}_BP$  is not
necessarily a Lie group, nor does it satisfy the conditions in \cite{LU14}, and we have not 
developed the classification theory in such full generality.

In \autoref{sec:free-loop-spaces}, we show that there is a weak $G$-equivalence between the free loop space $LB_G\Pi $ and the
adjoint bundle $ Ad(E_G\Pi) : =E_G \Pi \times_{\Pi} \Pi_{\mathrm{ad}}$ as $G$-fibrations over
$B_G\Pi$ (\autoref{thm:freeLoop}).

\textbf{Organization of the paper.} We give preliminaries of $G$-CW
complexes in \autoref{sec:g-cw-complexes}.  We review non-equivariant bundles in
\autoref{sec:non-equiv-bundl}. We give definitions of equivariant bundles and
compare the definitions in
\autoref{sec:defin-equiv-bundl} - \autoref{subsec:split-ext}. We study the
fixed points of equivariant bundles in \autoref{sec:fixed-point-theorems}. We
prove the aforementioned theorems about classifying spaces in \autoref{sec:classify}.
 
\textbf{Acknowledgements.}
This is part of the author's PhD thesis. The author is indebted to her advisor, Peter May, for explaining
his past work and for his enormous help with the writing. The
author thanks the referee for the useful comments and writing suggestions.

\section{Equivariant bundles}
\label{sec:equiBundle}
\subsection{$G$-CW complexes}
\label{sec:g-cw-complexes}
In this subsection, we give some preliminaries.
For a compact Lie group $G$, a $G$-CW complex $X$ is a union of $G$-spaces
$X^n$ constructed as follows. We start with a disjoint union of orbits $X^0 =
\sqcup_{K} G/K $. Each $X^n$ is inductively obtained by gluing $n$-cells to
$X^{n-1}$.  An $n$-cell is of the form
$G/K \times D^n$  where $K \subgroup G$ is  a closed subgroup and $D^n = \{
\mathbf{x} \in \bR^n| | \mathbf{x} | \leq 1 \}$ is a disk. There are $G$-maps
from $\sqcup_{K} G/K \times S^{n-1}$, the boundary of the $n$-cells, to $X^{n-1}$ and
$X^n = X^{n-1} \cup_{(\sqcup_{K} G/K \times S^{n-1})} \big(\sqcup_{K} G/K \times D^n\big)$.

\begin{defn}
  A map $f: X \to Y$ between spaces is said to be a weak equivalence if it
  induces a bijection $\pi_0(X) \to \pi_0(Y)$ and an isomorphism $\pi_n(X,x) \to
  \pi_n(Y, f(x))$ for any $x\in X$ and $n \geq 1$.
    A $G$-map  $f:X \to Y$  between $G$-spaces is said to be a weak $G$-equivalence
  if $X^H \to Y^H$ is a weak equivalence for any subgroup $H \subset G$.
\end{defn}

A $G$-space $X$ is said to have the homotopy type of a $G$-CW complex if there
is a $G$-CW complex $X'$ and a $G$-homotopy equivalence $X' \simeq X$.
\begin{thm}[Equivariant Whitehead theorem]
  A weak $G$-equivalence between $G$-spaces having homotopy types of $G$-CW
  complexes is a $G$-homotopy equivalence.
\end{thm}
This theorem allows us to use the induced maps on homotopy groups to detect
homotopy equivalences when working with $G$-CW complexes.

\subsection{Non-equivariant bundles}
\label{sec:non-equiv-bundl}
We start with a review of non-equivariant bundles.

  A \emph{fiber bundle} with fiber $F$  is a map $p:E \to B$ with an open cover
  $\{U_i\}$ of $B$ and homeomorphisms $\phi_i: p^{-1}(U_i) \cong U_i \times F$. The $U_i$ are called
  \emph{coordinate neighborhoods} and the $\phi_i$ are called \emph{local trivializations}.

  The structure group of a fiber bundle gives information about the change of local trivializations 
  under changes of coordinate neighborhoods.
  Let $\Pi$ be a topological group with an effective action on $F$.
  Here, \emph{effective} means $\Pi \to \mathrm{Aut}(F)$ is an injection.
  A bundle with fiber $F$ is said to have \emph{structure group} $\Pi$, if  
  for any two coordinate neighborhoods with $U_i \cap U_j \neq \varnothing$, the
  composite $\phi_i\phi_j^{-1}: (U_i \cap U_j ) \times F \to (U_i \cap U_j
  ) \times F$ is given by $(b,f) \mapsto (b,g_{ij}(b)(f))$ for some continuous
  function $g_{ij}: U_i \cap U_j \to \Pi$, called a \emph{coordinate transformation}.
  We always topologize $\mathrm{Aut}(F)$ with the compact-open topology of mapping spaces.
    If $F$ is a compact Hausdorff space, $\mathrm{Aut}(F)$ is a topological group; if $F$ is only
    locally compact, there are more technical assumptions for the inverse map to be continuous due to Arens
    (See \cite[I.5.4]{Steenrod}).  Morally, a
  fiber bundle with fiber $F$ is automatically a fiber bundle with the implicit structure
  group $\mathrm{Aut}(F)$. Having an explicit structure group $\Pi$ is extra data
  to reduce  the structure group to a smaller one.

  One can associate a principal $\Pi$-bundle to a fiber bundle with structure group $\Pi$.
  An \emph{admissible map} of the bundle is a homeomorphism
  $\psi: F \to  p^{-1}(b)$ for some $b \in U_i$, satisfying $\phi_i\psi \in \Pi$.
  The set of admissible maps has a natural topology, and the \emph{associated principal $\Pi$-bundle} of $p$ is the space of admissible maps.
  
  The following immediate observation about admissible maps hides the local trivializations in the
  background.
  \begin{lem}
    \label{lem:admissble}
     A map $\psi: F \to F_b$ is admissible if and only if for any admissible
     map $\zeta : F \to F_b$, the composite $\zeta^{-1}\psi$ is in~$\Pi$. \qed
   \end{lem}

  Let ${p_1:E_1 \to B_1}$ and ${p_2: E_2 \to B_2}$ be two fiber bundles with fiber $F$ and structure
  group $\Pi$.
  A \emph{morphism} between them is a bundle map $\chi: E_1 \to E_2$ such that for any local trivializations $\phi_U:
  p_1^{-1}(U) \cong U \times F$ and $\phi_V: p_2^{-1}(V) \cong V \times F$, the composite
\begin{equation}
\label{eq:bundle-morphism}
\phi_V \chi  \phi_U^{-1}: (U \cap \chi^{-1}(V) ) \times F \to (\chi(U) \cap V ) \times F
\end{equation}
  is given by $(b,f) \mapsto (\chi(b),g_{VU}(b)(f))$, where $g_{VU}: U \cap \chi^{-1}(V) \to \Pi$ is
  some continuous function.
  Such a morphism induces a morphism between the two associated principal $\Pi$-bundles.

  \begin{ass}
    \label{assumption:Steenrod}
    We always assume that $\Pi$ has the subspace topology of $\mathrm{Aut}(F)$.
  \end{ass}
  We pause to explain the role of \autoref{assumption:Steenrod}. Suppose $\chi$
  is a bundle map and we would like to check if $\chi$ is 
  a morphism, that is, whether it respects the structure group. It seems as if one only needs to check that $\chi$ sends
  all admissible maps to admissible maps.
  This is not true without \autoref{assumption:Steenrod}, since the set of admissible maps does not see the
  topology.  Steenrod \cite[I.5]{Steenrod}
  studied this difference carefully and concluded that
  \autoref{assumption:Steenrod} will
  resolve the discrepancy.

 We include some explanation here for completeness:
  What the set of admissible maps sees is an Ehresmann-Feldbau bundle with structure
  group $\Pi$, which has now become an obsolete notion.
  An Ehresmann-Feldbau bundle is a bundle $p:E \to B$ with fiber $F$ and a set of homeomorphism
  $\psi: F \cong p^{-1}(b)$ for all $b \in  B$, called admissible maps.
  It is required that for any $b \in U_i$, the composite ${F = \{b\} \times F \to U_i \times F
  \overset{\phi_i^{-1}}{ \to } p^{-1}(U_i)} $ is admissible, and that for any $b \in B$ and
  any admissible map $\psi: F \to p^{-1}(b)$, all the admissible maps $F \to
  p^{-1}(b)$ are exactly $\psi \circ \ele$ for some $\ele \in \Pi$.
  {While this aligns with \autoref{lem:admissble} when the bundle has a structure group $\Pi$,
  there is a difference between the two notions, which lies exactly in that an Ehresmann-Feldbau
  bundle does not require $\Pi$ to have a topology. In other words, the coordinate transformations
    $g_{ij}$ are not asked to be continuous, which
  is equivalent to putting the trivial topology on $\Pi$.
  If $\Pi$ does start life with a different topology, the coordinate transformations
  $g_{ij}$ obtained from an Ehresmann-Feldbau bundle may not be continuous.
  It is shown in \cite[I.5.4]{Steenrod} that with
  \autoref{assumption:Steenrod}, the $g_{ij}$'s are automatically
  continuous, so that a fiber bundle with local trivializations has structure group $\Pi$ if and only if
  the admissible maps satisfy \autoref{lem:admissble}.
  We have the following criteria:
   \begin{prop}
     \label{rem:admissible-criteria} 
  A bundle map $\chi: E_1 \to E_2$ is a morphism of fiber bundles with
  structure group $\Pi$ if and only if either of the two equivalent conditions is true:
\begin{enumerate}
\item \label{item:admissible1} If $F_1$ is a fiber in $E_1$ and $F_2$ is a fiber in $E_2$ such that
  $\chi$  maps $F_1$ to $F_2$, then the composite  $\zeta^{-1} \chi \psi$ is in $\Pi$ for any
  admissible maps $\psi: F \to F_1$ and $\zeta: F \to F_2$.
\item \label{item:admissible2}For any admissible map $\psi: F \to F_1$ to a fiber in $E_1$,
  the composite $\chi \psi$ is an admissible map to a fiber in $E_2$.
\end{enumerate}
   \end{prop}
   \begin{proof}
     We need to check that for any $\phi_U$, $\phi_V$ as in \autoref{eq:bundle-morphism}, the
     desired $g_{VU}$ exists. With \autoref{assumption:Steenrod}, it suffices to check that for any
     $b \in U \cap \chi^{-1}(V)$, there exists a desired $g_{VU}(b) \in \Pi$.
     This is part~\autoref{item:admissible1}. Part~\autoref{item:admissible2} follows from
     \autoref{lem:admissble}.
   \end{proof}

   \begin{exmp}
    The most familiar case is when $F$ is a vector space ($\bR^n$ or $\bC^n$)
    and $\Pi=GL_n$ is the corresponding general linear group.
    By definition of the general linear group, $\chi$ being a bundle map is
    equivalent to it being fiberwise linear and non-degenerate.
  \end{exmp}
   
  The following well-known structure theorem turns the problem of classifying fiber
  bundles into classifying principal bundles.
  \begin{thm}
    \label{thm:structure}
    Let $\Pi$ be a compact Lie group. Let $B, F$ be spaces. Assume that $\Pi$ acts effectively on $F$.
    Then there is an equivalence of categories between
    \{fiber bundles over $B$ with fiber $F$ and structure group $\Pi$\}
    and \{principal $\Pi$-bundles over $B$\}.
  \end{thm}
  \begin{proof}
 We have already shown how to construct a principal $\Pi$ bundle from a fiber bundle with fiber
     $F$ and structural group $\Pi$ at the beginning of this subsection.
  In the other direction, given a principal $\Pi$-bundle $P \to B$,
  the map $P\times_{\Pi} F \to B$ is a fiber bundle with fiber $F$ and structure group
  $\Pi$. These two constructions are functorial and inverse of each other.
  Indeed, \cite[I]{Steenrod} described both types of bundles using local transformations, called
  coordinate bundles, where the equivalence becomes transparent.
  \end{proof}

  \subsection{Definitions of equivariant bundles}
  \label{sec:defin-equiv-bundl}
  When it comes to the equivariant story, there are definitions of
  different generality, both on the fiber bundle side and on the principal
  bundle side. The reason is that the ambient group $G$ could interact
  non-trivially with the structure group $\Pi$.
  We start with the simplest definition where ``$G$ and $\Pi$ commute'' 
  \cite{L82}.
  
\begin{defn}
  \label{defn:Gvector}
  A $G$-fiber bundle with fiber $F$ and structure
  group $\Pi$ is a map $p:E \to B$ such that the following statements hold:
  \begin{enumerate}
    \item The map $p$ is a non-equivariant fiber bundle with fiber $F$ and
      structure group~$\Pi$;
    \item Both $E$ and $B$ are $G$-spaces and $p$ is $G$-equivariant;
    \item \label{item:Gvector3} The $G$-action is given by morphisms of bundles with structure group $\Pi$.
    \end{enumerate}
\end{defn}

\begin{prop}
  \label{rem:Gvector-admissible}
  The requirement in \autoref{item:Gvector3} above is equivalent to the following: for any $g \in G$
  and admissible map $\psi: F \to F_b$, the composite
  $F \overset{\psi}{ \to } F_b \overset{g}{ \to } F_{gb}$ is also admissible.
\end{prop}
\begin{proof}
  By \autoref{rem:admissible-criteria}.
\end{proof}

  \begin{rem}
    Let $G$ be a finite group. We take $F=\bR^n$ and $\Pi = \mathrm{GL}_n(\bR)$  in
    \autoref{defn:Gvector}.  Although $\mathrm{GL}_n(\bR)$ is not compact, the definition still
    works, and we obtain a $G$-$n$-vector bundle.
  \end{rem}

\begin{defn}
  \label{defn:Gprincipal}
  A principal $G$-$\Pi$-bundle is a map ${p:P \to B}$ such that the following statements hold:
    \begin{enumerate}
      \item The map $p$ is a non-equivariant principal $\Pi$-bundle;
      \item Both $P$ and $B$ are $G$-spaces and $p$ is $G$-equivariant;
      \item The actions of $G$ and $\Pi$ commute on $P$.
    \end{enumerate}
\end{defn}

\begin{rem}
  This is called a principal $(G,\Pi)$-bundle in \cite[IV1]{LMS86}.
\end{rem}
 As in the non-equivariant case, we write the $\Pi$-action on the right of a
  principal $G$-$\Pi$-bundle $P$; for convenience of diagonal action, we consider
  $P$ to have a left $\Pi$-action, that is, $\ele \in \Pi$ acts on $z \in P$ by
  $\ele z = z\ele^{-1}$.

The structure theorem formally passes to this equivariant context.
  \begin{thm}
\label{thm:G-structure-1}
    Let $G,\Pi$ be compact Lie groups and $F,B$ be spaces. Assume that $\Pi$ acts effectively on $F$.
    Then there is an equivalence of categories between
    \{$G$-fiber bundles over $B$ with fiber $F$ and structure group $\Pi$\}
    and \{principal $G$-$\Pi$-bundles over $B$\}.
  \end{thm}
  \begin{proof}
    The two types of $G$-bundles in Definitions \ref{defn:Gvector} and \ref{defn:Gprincipal}
    are indeed objects with a
    $G$-action in the corresponding non-equivariant category.
    So  the equivalence in the non-equivariant structure theorem
    restricts to give an equivalence on the $G$-objects.
  \end{proof}

  \medskip
However, Definitions \ref{defn:Gvector} and \ref{defn:Gprincipal} are not ideal for studying some interesting
cases. In the most general scenario, we want to study a map $p:E \to B$ that happens
to be both a fiber bundle with structure group $\Pi$ and a $G$-map between
$G$-spaces. It is true that $p$ is a $G$-fiber bundle with structure group
$\mathrm{Aut}(F)$, but $p$ is usually not a $G$-fiber bundle with structure group $\Pi$.
In other words, we can not reduce the structure group even though we know non-equivariantly it reduces to $\Pi$.
Below, we give two concrete examples of this sort.

The first example is Atiyah's Real vector bundles \cite{Atiyah}.
\begin{exmp} \label{eg:Real}
Let $G=C_2$. A Real vector bundle is a map $p:E \to B$ such that 
\begin{itemize}
\item The map $p$ is a complex vector bundle of dimension $n$;
\item The non-trivial element of $C_2$ acts anti-complex-linearly.
\end{itemize}
In this case, $p$ is a $C_2$-bundle with structure group $O(2n)$, but not $U(n)$.
\end{exmp}

The second simple, but illuminating, example is from \cite{LMS86}.
\begin{exmp} \label{eg:Gbundle}
For $G$-spaces $B$ and $F$, the projection $p: B \times F \to B$ is not a $G$-bundle with structure
group $e$ unless $G$ acts trivially on $F$.  
\end{exmp}
\begin{proof}
  The admissible maps for $p$ are only the inclusions of fibers
  $${\psi_b: \{b\} \times F \to B \times F}.$$
  An element $g \in G$ acts by a bundle map if and
  only if for all $b$, the composite 
\begin{equation*}
  \{b\} \times F  \overset{\psi_b}{ \to }  p^{-1}(b)
  \overset{g}{ \to } p^{-1}(gb) \overset{\psi_{gb}^{-1}}{ \to } \{gb\}  \times F
\end{equation*}
  is  in the structure group. But this map is merely the $g$ action on $F$.
\end{proof}

\medskip
Consequently, we would like a more general version than Definitions \ref{defn:Gvector} and
\ref{defn:Gprincipal}.
To work with Real vector bundles, Tom~Dieck \cite{TD} introduced a complex conjugation action of $C_2$
  on $U(n)$. Lashof--May \cite{LM86} had the idea to further introduce a total group that is
  the extension of the structure group $\Pi$ by $G$.
Tom~Dieck's work became a special case of a split extension, or equivalently a semidirect product.
One good, but rather brief and sketchy, early reference for both is \cite[IV1]{LMS86}; we shall
  flesh out that source and come back to the two examples afterwards.

We start with the well-studied principal bundle story.

\begin{defn}(\cite{LM86})
\label{defn:GGammaprincipal}
  Let $1 \to \Pi \to \Gamma \to G \to 1$ be an extension of compact Lie groups. 
  A principal $(\Pi;\Gamma)$-bundle is a map $p:P \to B$ such that the
  following statements hold:
    \begin{enumerate}
      \item The map $p$ is a non-equivariant principal $\Pi$-bundle;
      \item The space $P$ is a $\Gamma$-space; $B$ is a $G$-space. Viewing $B$ as a
 $\Gamma$-space by pulling back the action, the map $p$ is $\Gamma$-equivariant.
    \end{enumerate}
\end{defn}

\begin{rem}
  The total space $P$ does not have a $G$-action in general.
  It only does so when we specify a splitting $G \to \Gamma$.
  An example of this sort will be discussed in \autoref{subsec:split-ext}.  
\end{rem}

  \begin{defn}
  A morphism between two principal $(\Pi;\Gamma)$-bundles $p_1: P_1 \to B_1$ and $p_2: P_2 \to
    B_2$ is a pair of maps $(\bar{f},f)$ fitting in the commutative diagram  
\begin{equation*}
  \begin{tikzcd}
    P_1 \ar[r,"{\bar{f}}"] \ar[d,"{p_1}"'] & P_2 \ar[d,"{p_2}"] \\
    B_1 \ar[r,"f"] & B_2
  \end{tikzcd}
\end{equation*}
such that $f$ is $G$-equivariant and $\bar{f}$ is $\Gamma$-equivariant.
  \end{defn}

  \begin{exmp}
    Let $y \in \Gamma$ be with image $g \in G$. The action map $(y,g)$ is an automorphism.
  \end{exmp}

Taking $\Gamma = \Pi \times G$, we recover the principal $G$-$\Pi$-bundles
of \autoref{defn:Gprincipal}.
In this case we have two names for the same thing. This could be confusing,
but since a ``principal $G$-$\Pi$-bundle'' looks more natural than a ``principal
$(\Pi; \Pi \times G)$-bundle'' for this thing, we will keep both names.

 Taking $\Gamma$ to be a split extension, or equivalently, $ \Gamma = \Pi \rtimes_{\alpha} G$ for some
 group homomorphism $\alpha:G \to \mathrm{Aut}(\Pi)$,
 we recover Tom~Dieck's principal $(G,\alpha,\Pi)$-bundles.

 \begin{rem}
   \label{rem:semidirect-group}
 To be useful later, we write the elements of $\Gamma = \Pi \rtimes_{\alpha} G$ as $(\ele,g)$
 for $\ele \in \Pi, g \in G$ and write $\alpha(g) \in \mathrm{Aut}(\Pi)$ as $\alpha_g$.
 We have the following facts:
 \begin{itemize}
   \item The product in $\Gamma$ is given by $(\ele,g)(\elm,h) = (\ele\alpha_{g}(\elm), gh)$ 
 (That is, $g$ acts on $\elm$ when they interchange);
   \item The identity element is $(e,e)$;
   \item The inverse is $(\ele,g)^{-1} = (\alpha_{g^{-1}}(\ele^{-1}), g^{-1} )$;
   \item The elements $(e,g)$ form a subgroup of $\Gamma$ that is canonically
     isomorphic to $G$;
   \item A space with $\Gamma$-action is a space with both $\Pi$ and $G$ actions
     such that  $$g(\ele(-))=\alpha_{g}(\ele)(g(-)), \text{ which is indeed }(\alpha_{g}(\ele),g) (-).$$
 \end{itemize}
\end{rem}

 The fiber bundle story is not as clear. It turns out that 
the appropriate fiber of an equivariant fiber bundle is not just the preimage of
any point, but rather with a preassigned action of $\Gamma$.
This is unnatural at first glance, for example in a $G$-vector bundle we will not
 expect there to be an $(O(n) \times G)$-action on the fiber $\bR^n$.
 We will explain why this is necessary at the end of this subsection.
 How $G$-vector bundles fit in this context will be stated in \autoref{ex:G-vector-bundle}.
 Let us start with the definition:
 \begin{defn}(\cite[IV1]{LMS86})
   \label{defn:Gfiber}
   Let $1 \to \Pi \to \Gamma \to G \to 1$ be an extension of compact Lie groups
   and $F$ be a space with $\Gamma$-action.
   A $G$-fiber bundle with fiber $F$, structure group $\Pi$ and total group $\Gamma$
   is a map $p:E \to B$ such that the following statements hold:
  \begin{enumerate}
    \item The map $p$ is a non-equivariant fiber bundle with fiber $F$ and structure group~$\Pi$;
    \item Both $E,B$ are $G$-spaces and $p$ is a $G$-map;
    \item \label{item:Gfiber-admissible}
      For any  $g \in G$ and admissible maps
      $\psi: F \to F_b$ and $\zeta: F \to F_{gb}$, the composite      
      \begin{equation*}
        F \overset{\psi}{ \to } F_b \overset{g}{ \to } F_{gb} \overset{\zeta^{-1}}{ \to } F
      \end{equation*}
      is a lift $y \in \Gamma$ of $g \in G$. In other words, the $y$ in the following
      diagram is asked to be a lift of $g \in G$ in $\Gamma$:
\begin{equation*}
 \begin{tikzcd}
        F \ar[d,"\psi"',"\cong"] \ar[r,dotted, "y"] & F \ar[d, "\zeta","\cong"'] \\
        F_b \ar[r,"g"]& F_{gb}
      \end{tikzcd} 
\end{equation*}
\end{enumerate}
\end{defn}

\begin{prop}
  \label{rem:Gfiber-admissible}
  The requirement \autoref{item:Gfiber-admissible} above is equivalent to the following:
      For each $y \in \Gamma$ with image $g \in G$ and admissible map
      $\psi: F \to F_b$, the composite      
\begin{equation*}
F \overset{y^{-1}}{ \to } F \overset{\psi}{ \to } F_b \overset{g}{ \to } F_{gb}
\end{equation*}
is also admissible.
\end{prop}
\begin{proof}
  For any  two lifts $y$ and $y'$ of $g$, $y'y^{-1}$ is a lift of $e \in G$, so it is in $\Pi$.
The claim then follows from \autoref{lem:admissble}.
\end{proof}

Taking $g=e$ in \autoref{rem:Gfiber-admissible}, the possible lifts $y$ are exactly the elements of $\Pi$,
so we just see the non-equivariant structure group (compare with
  \autoref{lem:admissble}); taking general $g$, the assignment $\psi \mapsto g \psi y^{-1}$ is
  mimicking the action by an element of $\Pi$ on  the admissible map $\psi$, but it changes the
  fiber from over $b$ to over $gb$.
  In this sense, the extension of  the structure group $\Pi$ to  the
  total group $\Gamma$ is used to regulate admissible maps to fibers over the orbit of $b$.

 \begin{defn} \label{defn:morphismFiber}
Let  $p_1: E_1 \to B_1$ and $p_2: E_2 \to  B_2$ be two $G$-fiber bundles with fiber $F$,
  structure group $\Pi$ and total group $\Gamma$.  A morphism between them is a pair of maps
  $(\bar{f},f)$ fitting in the commutative diagram 
\begin{equation*}
  \begin{tikzcd}
    E_1 \ar[r,"{\bar{f}}"] \ar[d,"{p_1}"'] & E_2 \ar[d,"{p_2}"] \\
    B_1 \ar[r,"f"] & B_2
  \end{tikzcd}
\end{equation*}
such that the following statements hold:
  \begin{enumerate}
  \item \label{item:morphism1} The pair $(\bar{f},f)$ is a non-equivariant morphism between bundles with fiber $F$ and structure group $\Pi$.
  \item Both $\bar{f}$ and $f$ are $G$-equivariant.
    \end{enumerate}
  \end{defn}
\begin{rem}
  By \autoref{rem:admissible-criteria}, the condition \autoref{item:morphism1}
  of \autoref{defn:morphismFiber}
  is explicitly the following: For any admissible map $\psi: F \to F_1$ to a fiber in $E_1$,
  the composite $\bar{f} \psi$ is an admissible map to a fiber in $E_2$.
\end{rem}
We do not have a
  requirement on a morphism regarding the condition~\autoref{item:Gfiber-admissible} of
  \autoref{defn:Gfiber} because it is automatic: if $\psi$ is
  admissible, we have that $g\psi y^{-1}$ is admissible and so is $\bar{f}(g \psi y^{-1})$. But
   $\bar{f}g = g \bar{f}$, so $g (\bar{f} \psi) y^{-1}$ is also admissible.

   \medskip
   As opposed to \autoref{defn:Gvector}, in \autoref{defn:Gfiber} the $\Gamma$-action on the total
   space $E$ can restrict to
   a $G$-action only when there is a splitting of the extension given by $G \to \Gamma$.
   The following example illustrates that varying the splitting map can give different  $G$-fiber
   bundle descriptions of the same bundle. It will be discussed in \autoref{subsec:split-ext}.
\begin{exmp}
     \label{ex:G-vector-bundle}
      A $G$-$n$-vector bundle is both a $G$-fiber bundle with fiber $\bR^n$, structure group
  $O(n)$ and total group $O(n) \times G$ and a $G$-fiber bundle with fiber $V$, structure group
  $O(V)$ and total group $O(V) \rtimes G$. (Here, we take $\Gamma = O(n) \times G \cong O(V) \rtimes G$.)
\end{exmp}

\begin{exmp}
  \label{example:real}
  A Real vector bundle is a $C_2$-fiber bundle with fiber $\bC^n$, structure
  group $U(n)$ and total group $\Gamma = U(n) \rtimes_{\alpha} C_2$, where $\alpha: C_2 \to
  \mathrm{Aut}(U(n))$ sends the non-trivial element of $C_2$ to the entry-wise complex-conjugation of $U(n)$.
\end{exmp}
\begin{proof}
  Let the non-trivial element $a$ of $C_2$ act by complex conjugation on $\bC^n$. This extends the
  $U(n)$-action to a $\Gamma$-action by \autoref{rem:semidirect-group}. 
We only need to check that \autoref{defn:Gfiber}~\autoref{item:Gfiber-admissible} holds for $g=a$. An automorphism $X$ of $\bC^n$ is
anti-complex-linear if and only if $A = X\circ a$, the pre-composition of $X$ with conjugation,
is complex-linear. So $A$ is an element of $U(n)$, and $X=(A,a)$ is the lift of $a$ in $U(n) \rtimes_{\alpha} C_2$.
\end{proof}

\begin{exmp}
  For $G$-spaces $B$ and $F$, the projection $B \times F \to B$
  is a $G$-fiber bundle with fiber $F$, structure group $e$ and total group $\Gamma = G$.
\end{exmp}
\begin{proof}
  The proof in \autoref{eg:Gbundle} verifies \autoref{defn:Gfiber}~\autoref{item:Gfiber-admissible}.
\end{proof}

When $\Gamma = \Pi \times G$, \autoref{defn:Gvector} can be considered as a
special case of \autoref{defn:Gfiber}: The canonical way to extend the
$\Pi$-action on $F$ to a $\Gamma$-action is by taking the trivial
$G$-action. With this convention, a $G$-fiber bundle in the first definition
satisfies the second definition, as we will show shortly in \autoref{prop:bundle-action1}. In fact, there could be multiple ways to
extend the action, such as in \autoref{ex:G-vector-bundle}.
On the other hand, \autoref{eg:Gbundle} shows that even for a trivial bundle $B
\times F \to F$, the more general \autoref{defn:Gfiber} is needed if $G$ acts non-trivially on $F$.

We have the following structure theorem in the context of Definitions
\ref{defn:GGammaprincipal} and \ref{defn:Gfiber}:

\begin{thm}(\cite[IV1]{LMS86})
\label{thm:G-structure-2}
  For any $\Pi$-effective $\Gamma$-space $F$ and $G$-space $B$,
  there is an equivalence of categories between \{$G$-fiber bundles with structure group $\Pi$, total group $\Gamma$ and
  fiber $F$ over $B$\} and \{principal $(\Pi;\Gamma)$-bundles over $B$\}. 
\end{thm}
\begin{proof} This is an expansion of the sketchy proof in the reference. For brevity, we
    refer to the two categories as equivariant fiber bundles and equivariant principal bundles.
  
  Given an equivariant fiber bundle $E \to B$, we take the non-equivariant associated principal bundle $\mathrm{Fr}_F(E) \to B$.
  It suffices to give a $\Gamma$-action on $\mathrm{Fr}_F(E)$ such that $\mathrm{Fr}_F(E) \to B$ is a $G$-map.
  For $y \in \Gamma$ with image $g \in G$ and an admissible map $\psi : F \to F_b$,
  let $y (\psi) = g \psi y^{-1}$. By \autoref{rem:Gfiber-admissible}, $g \psi y^{-1}$ is an admissible map to
  the fiber over $gb$. This shows that $\mathrm{Fr}_F(E) \to B$ is an equivariant principal bundle.

  Given an equivariant principal bundle $P \to B$, let $E = (P \times F)/\Pi \to B$ be the fiber bundle with
  admissible maps $\psi_p: F \to E$ of the form $\psi_p(f)=[p,f]$ for some $p \in P$.
  We verify the three conditions for $E \to B$ to be an equivariant fiber bundle.
  Firstly, $E \to B$ is a non-equivariant fiber bundle with structure group $\Pi$.
  Secondly, we describe the $G$-action
  on $E$. Take the diagonal $\Gamma$-action on $P \times F$.
  For any space with $\Gamma$-action $X$, we can define a $\Gamma/\Pi \cong G$-action on 
   $X/\Pi$ by lifting $g \in G$ to $y \in \Gamma$ and let $g[x] = [yx]$ for $x \in X$.
  Since $\Pi$ is a normal subgroup of $\Gamma$, this is a well-defined action, independent of
  choice of $y$ or representative $x$.
  For $X = P \times F$, this gives $(P \times F)/\Pi$ a $G$-action. Since $P \to B$ is
  $\Gamma$-equivariant, it can be checked that $E \to B$ is $G$-equivariant.
  Thirdly, we show that the condition in \autoref{rem:Gfiber-admissible} is satisfied.
  In fact, for $y \in \Gamma$ lifting $g \in G$ 
   and $p \in P$, we have  $g\psi_py^{-1} = \psi_{yp}$. To see this, evaluating on $f \in F$, we have
   \begin{align*}
     g\psi_py^{-1}(f) &  = g[p, y^{-1}f] & \text{ definition of $\psi$;} \\
& = [yp, yy^{-1}f] & \text{ definition of $G$-action;}\\
& = [yp,f] 
 = \psi_{yp}(f) & \text{ definition of $\psi$.}
\end{align*}

  These two constructions give inverse functors.  Given an equivariant fiber bundle $E \to B$, we have a map
\begin{equation*}
\xi: (\mathrm{Fr}_F(E) \times F)/\Pi \to E,\ \xi([\psi,f])= \psi(f).
\end{equation*}
Non-equivariantly we already know that $(\xi, \mathrm{id}_B)$ is a morphism of fiber bundles with
structure group $\Pi$ and that $\xi$ is a homeomorphism.
To check that $\xi$ is $G$-equivariant, suppose $g \in G$ lifts to $y \in \Gamma$. Then
\begin{equation*}
  g([\psi, f]) = [y(\psi),yf] = [g\psi y^{-1},yf]
\end{equation*}
and $\xi([g\psi y^{-1},yf]) = (g \psi y^{-1})(yf) = g(\psi(f))$.
So $(\xi, \mathrm{id}_B)$ is a morphism of equivariant fiber bundles
by \autoref{defn:morphismFiber}. It is an isomorphism because the non-equivariant inverse is also an
equivariant inverse as it is a homeomorphism.
Given an equivariant principal bundle $P \to B$, we have a map which we abusively denote by
\begin{equation*}
\psi: P \to \mathrm{Fr}_F((P \times F)/\Pi), \ p \mapsto \psi_p.
\end{equation*}
Here, $\psi_p$ is the previously defined admissible map of $(P \times F)/\Pi$, thus an element of
its associated principal bundle.  Again, non-equivariantly we know that the map $\psi$ is a homeomorphism
(the $\Pi$-effectiveness is needed to assure that if $p \neq q$ in $P$, then $\psi_p \neq \psi_q$).
To check that $\psi$ is $\Gamma$-equivariant,
the definition of the $\Gamma$-action on admissible maps gives $y\psi_p = g\psi_py^{-1}$ and we have
verified $g\psi_py^{-1} = \psi_{yp}$, so we have $y\psi_p = \psi_{yp}$. Thus, $(\psi, \mathrm{id}_B)$ is a
morphism of equivariant principal bundles. It is also an isomorphism.
\end{proof}

\begin{rem}
  The isomorphisms $\xi$ and $\psi$ in the proof are natural and provide the unit and counit maps of the
  adjunction 
\begin{equation*}
\mathrm{Hom}((P \times F)/\Pi, E) \cong \mathrm{Hom}(P, \mathrm{Fr}_F(E))
\end{equation*}
\begin{equation*}\small
  \begin{tikzcd}
    (- \times F)/\Pi:
    \left\{ \text{
      \begin{tabular}[h]{c}
        principal $(\Pi;\Gamma)$- \\
        bundles over  $B$
      \end{tabular} } \right\}
  \ar[r, shift left] &
    \left\{ \text{
      \begin{tabular}[h]{c}
        $G$-fiber bundles  over $B$ \\
        with structure group $\Pi$, \\
        total group $\Gamma$ 
        and fiber $F$
      \end{tabular} } \right\}  :
    \mathrm{Fr}_F(-)  \ar[l, shift left ]
  \end{tikzcd}
\end{equation*}
\end{rem}

We can see in the proof  of \autoref{thm:G-structure-2}
that it is essential for $F$ to have a $\Gamma$-action.
If $P$ were a principal $(\Pi;\Gamma)$-bundle and the fiber $F$ only had a $\Pi$-action,
the associated fiber bundle $(P \times F)/\Pi$ would not have a $G$-action.
If we insist on our notion of a $G$-fiber bundle to be a $G$-map between $G$-spaces, this is the
price to pay.

\subsection{Comparisons of definitions}
\label{sec:comp-definition}
We have two concepts of $G$-fiber bundles. One is the $G$-fiber bundle with fiber $F$ and structure
group $\Pi$ as in \autoref{defn:Gvector};
the other is the $G$-fiber bundle with fiber $F$, structure group $\Pi$ and total group $\Gamma$ for
a specific extension of compact Lie groups  $1 \to \Pi \to \Gamma \to G \to 1$, as in
\autoref{defn:Gfiber}.
The differences between the concepts are two-fold: in the first one, $G$ acts by bundle maps, but in
the second one, the $G$-action is regulated by $\Gamma$;
in the first one, $F$ has only a $\Pi$-action, but in the second one, $F$ has a
$\Gamma$-action. We compare these two concepts and show that
the first concept is a special case of the second when  $\Gamma \cong \Pi \times G$ and 
$\Gamma$ acts on $F$ via the projection $\Pi \times G \to \Pi$
(\autoref{prop:bundle-action1}).

We start with some simple group theory observations that will come into play.
\begin{defn}
  A retraction $\Gamma \to \Pi$ is a group homomorphism that restricts to  identity on the subgroup $\Pi$.
\end{defn}

It turns out that $\Gamma$ admits a retraction to $\Pi$ if and only if it is isomorphic
to $\Pi \times G$. We prove this explicitly in the case of a semidirect product first, then for general $\Gamma$.
\begin{prop}
  \label{ex:action-on-F-2}
  Let $\Gamma = \Pi \rtimes_{\alpha} G$ be a split extension.
  Then 
\begin{enumerate}
\item  \label{item:split1}The retractions $\tilde{\beta}: \Gamma \to \Pi$ are in bijection to
  homomorphisms $\beta: G \to \Pi$ satisfying $ \alpha_g(\ele) = \beta(g)\ele \beta(g)^{-1}$
  for all $g \in  G$ and $\ele \in \Pi$.
  (Note that for a given $\alpha: G \to \mathrm{Aut}(\Pi)$, the homomorphism $\beta$ may not exist.)
\item \label{item:split2} Each $\beta$ in \autoref{item:split1} specifies an isomorphism $ \Pi
  \rtimes_{\alpha} G \cong \Pi \times G$.
\end{enumerate} 
\end{prop}
\begin{proof}
  To see \autoref{item:split1}, we use the explicit expression for
  semidirect product group $\Gamma$ as in
  \autoref{rem:semidirect-group}. Suppose we have a retraction $\tilde{\beta}: \Gamma \to \Pi$.
  Let $\beta(g)$ be the image $\tilde{\beta}(e,g)$. Then $\beta$ is a group homomorphism. We
  have $\tilde{\beta}(\ele,e)=\ele$ and $$\tilde{\beta}(\ele,g) = \tilde{\beta}((\ele,e)(e,g)) =
  \ele\beta(g).$$
  In order for $\tilde{\beta}$ to be a homomorphism, it is required that the following two elements are equal
  for all $g,h \in G$ and $\ele,\elm \in \Pi$:   
  \begin{align*}
    \tilde{\beta}(\ele\alpha_g(\elm),gh) & = \ele\alpha_g(\elm)\beta(gh); \\
    \tilde{\beta}(\ele,g) \tilde{\beta}(\elm,h) & = \ele\beta(g)\elm\beta(h).
  \end{align*}
  Comparing the two lines gives $\alpha_g(\elm) = \beta(g)
  \elm\beta(g)^{-1}$. On the other hand, if we have $\beta$ as required, the
  formula $\tilde{\beta}(\ele,g) = \ele\beta(g)$ defines a retraction $\tilde{\beta}$.
  
  Given such a $\beta$, the group isomorphism in \autoref{item:split2} is given by
  \begin{equation*}
    \Pi \rtimes_{\alpha} G \cong \Pi \times G, \ (\ele,g) \mapsto (\ele\beta(g),g). \qedhere
  \end{equation*}
\end{proof}

\begin{prop}
  \label{prop:retraction}
  There is a  bijection of sets between $\{\text{retractions } \tilde{\beta}: \Gamma \to
  \Pi\}$ and $\{\text{isomorphisms of extensions } \Gamma \cong \Pi \times G \}$.
\end{prop}
\begin{proof}
   Consider $\Pi$ as a subgroup of
  $\Gamma$ and denote by $q$ the surjection $\Gamma \to G$. Given a retraction $\tilde{\beta}: \Gamma \to
  \Pi$, the map $(\tilde{\beta}, q): \Gamma \to \Pi
  \times G$ is a group isomorphism, and vice versa. 
  \begin{equation*}
  \begin{tikzcd}
  1 \ar[r] & \Pi \ar[d,equal] \ar[r] & \Gamma \ar[d, "{(\tilde{\beta}, q)}"] \ar[r, "q"]
    \ar[l,shift right, dotted, "{\tilde{\beta}}"']  & G \ar[d, equal] \ar[r] & 1 \\
  1 \ar[r] & \Pi \ar[r] &  \Pi \times G \ar[r] & G \ar[r] & 1
  \end{tikzcd} \qedhere
\end{equation*}
\end{proof}

We now compare Definitions \ref{defn:Gvector} and \ref{defn:Gfiber} in the following propositions.
Note that we can think about a retraction $\Gamma \to \Pi$ 
  as a chosen isomorphism $\Gamma \cong \Pi \times G$ of extensions by \autoref{prop:retraction}.
  \begin{prop}
    \label{prop:bundle-action1}
   Let $F$ be a space with an effective $\Pi$-action and $1 \to \Pi \to \Gamma \to G \to 1$ be an extension of
   compact Lie groups. Then one can extend the $\Pi$-action on $F$ to a $\Gamma$-action such that 
   a $G$-fiber bundle of \autoref{defn:Gvector} is always  a $G$-fiber bundle of
   \autoref{defn:Gfiber} if and only if there is a retraction $\Gamma \to \Pi$ and the
   extended $\Gamma$-action on $F$ is via the retraction.
  \end{prop}
  \begin{proof}
  Suppose we have $p: E \to B$ as in \autoref{defn:Gvector} and $F$ has an extended $\Gamma$-action.
  Then the only thing to check for $p$ to be a $G$-fiber bundle of \autoref{defn:Gfiber}
  is whether it satisfies the condition in
\autoref{rem:Gfiber-admissible}. That is, it suffices to
show for each $y \in \Gamma$ with image $g \in G$ and admissible homeomorphism
      $\psi: F \to F_b$, the composite $g\psi y^{-1}$ is also admissible.
      By \autoref{rem:Gvector-admissible},
$g \psi$ is admissible. So by \autoref{lem:admissble}, for $y \in \Gamma$, $g \psi y^{-1}$ is
admissible  if and only if $y$ acts on $F$ as an element in $\Pi$.
In other words, the group homomorphism $\Gamma \to \mathrm{Aut}(F)$ factors through $\Pi \to \mathrm{Aut}(F)$.
\end{proof}

The converse is also true.
\begin{prop}
  \label{prop:bundle-action2}
   Let  $1 \to \Pi \to \Gamma \to G \to 1$ be an extension of
   compact Lie groups and $F$ be a $\Pi$-effective $\Gamma$-space.
   Then a $G$-fiber bundle of \autoref{defn:Gfiber} is always
   a $G$-fiber bundle of \autoref{defn:Gvector} if and only if  $\Gamma$ acts on $F$ via a
   retraction $\Gamma \to \Pi$.
 \end{prop}
 \begin{proof}
   We can reverse the argument in \autoref{prop:bundle-action1}.
Suppose we have $p: E \to B$ as in \autoref{defn:Gfiber}; to check whether $p$ is a $G$-fiber bundle of
\autoref{defn:Gvector}, we only need to check whether the condition in \autoref{rem:Gvector-admissible} holds.
Take any admissible homeomorphism $\psi: F \to F_b$. By \autoref{rem:Gfiber-admissible},
for any $y \in \Gamma$ with image $g \in G$, $g \psi y^{-1}$ is admissible.
By \autoref{lem:admissble}, $g \psi$ is admissible if and only if $y$ acts on $F$ as an element in~$\Pi$.
 \end{proof}

Using Propositions \ref{prop:bundle-action1} and \ref{prop:bundle-action2}, we
can identity \autoref{defn:Gvector} as a special case of \autoref{defn:Gfiber}.
\begin{exmp}
  \label{ex:action-on-F-1}
  Let $\Gamma = \Pi \times G$ and $F$ be a space with an effective $\Pi$-action.
  We give $F$ the trivial $G$-action. Equivalently, this is viewing $F$
  as a space with $\Gamma$-action via the projection $\Gamma \to \Pi$.
  In this perspective, the structure theorem \autoref{thm:G-structure-1} is a special case of
  \autoref{thm:G-structure-2}.  
\end{exmp}
\begin{exmp}
  In particular, let $\Gamma = O(n) \times G$ and give $\bR^n$ the usual $O(n)$-action and
  the trivial $G$-action.
  We have an equivalence of the two concepts: 
\begin{itemize}
\item $G$-vector bundles with fiber $\bR^n$ (the classical $G$-equivariant vector bundles);
\item $G$-fiber bundles with
  fiber $\bR^n$, structure group $O(n)$ and total group $O(n)\times G$.
\end{itemize}  
 \end{exmp}

\begin{exmp}[non-example] \label{exmp:RealNot}
  For a Real vector bundle as in \autoref{example:real}, $\Gamma$ does not act on $\bC^{n}$ via
  $U(n)$ for any $n$.
  So a Real vector bundle is not a $C_2$-fiber bundle with fiber $\bC^n$ and structure group $U(n)$.
\end{exmp}
\begin{proof}
  There is no retraction $\Gamma \to U(n)$, because otherwise
  by \autoref{ex:action-on-F-2}, we would need an element $\beta(a)$ of $U(n)$
  such that $\beta(a) A = \bar{A} \beta(a)$
  for all $A \in U(n)$, where $\bar{A}$ is the complex conjugation of $A$.
  But this does not exist for any $n$.
\end{proof}

\subsection{Examples: the $V$-framing bundle}
\label{subsec:split-ext}
In the extension $1 \to \Pi \to \Gamma \to G \to 1$, the group $G$ is redundant because it is just
$\Gamma/\Pi$. However, due to the special role of the group $G$ in equivariant homotopy theory,
we would like to understand the $G$-action wherever applicable.
Since the total space of a principal $(\Pi; \Gamma)$-bundle
has only a $\Gamma$-action, we now focus on the case of split extensions, when we have a
specified group homomorphism $G \to \Gamma$.
This becomes relevant when we define and study the $V$-framing bundle of
a $G$-vector bundle for representations $V$. It turns out that $\mathrm{Fr}_V(E)$ and
$\mathrm{Fr}_{\bR^n}(E)$ are the same even as principal $(\Pi; \Gamma)$-bundles, but they have
different $G$-actions.

Let $F$ be a space with an effective $\Pi$-action and one can do some yoga with
the fiber. We fix a group homomorphism $\beta: G \to \Pi$. Then by \autoref{ex:action-on-F-2}, $\beta$ determines an isomorphism
\begin{equation}
  \label{eq:19}
\Pi \rtimes_{\alpha} G \cong \Pi \times G.
\end{equation} where $\alpha: G \to \mathrm{Aut}(\Pi)$ is the group
homomorphism given by 
\begin{equation}
\label{eq:20}
\alpha_g(\ele) = \beta(g) \ele \beta(g)^{-1}.
\end{equation}
We can let the groups in \autoref{eq:19} act on $F$ via their
retraction to $\Pi$. Note that this is the same abstract action but is
different for elements of the form $(\ele, g)$ on the two sides.
For clarity, we denote this space by $F'$. Explicitly,
$(\Pi \times G)$ acts on $F'$ by
$(\ele,g)(x) = \ele(x) \text{ for }x \in F'$;
$(\Pi \rtimes_{\alpha}G)$ acts on $F'$ by $
(\ele,g)(x) =  \ele \big(\beta(g)(x)\big) .
$
Inclusion to the second coordinate gives a canonical inclusion of $G$ into both
$\Pi \times G$ and $\Pi \rtimes_{\alpha}G$, but this is not compatible with the isomorphism
\autoref{eq:19}. The second image is the graph subgroup $\Lambda_{\beta} = \{(\beta(g), g) |
g \in G\} \subset \Pi \times G$. Consequently, the two $G$-actions on $F'$ are different.

In summary, we have an isomorphism of extensions in the situation, but it is not an isomorphism of
split extensions, as $(e,g)$ of $\Pi \rtimes_{\alpha} G$ is sent to $(\beta(g),g)$ in $\Pi \times G$.
\begin{equation*}
  \begin{tikzcd}
  1 \ar[r] & \Pi \ar[d,equal] \ar[r] & \Pi \rtimes_{\alpha} G \ar[d, "\cong ","\autoref{eq:19}"'] \ar[r]
  & G \ar[d,equal] \ar[r] \ar[l, shift right, dotted, "{(e,g) \mapsfrom g}"'] & 1 \\
  1 \ar[r] & \Pi \ar[r] & \Pi \times G \ar[r] & G \ar[r] \ar[l,  shift left,
  dotted,"{(e,g)\mapsfrom g}"]& 1
  \end{tikzcd}
\end{equation*}

As a consequence, we get the following trivial corollary of Propositions~\ref{prop:bundle-action1}
and \ref{prop:bundle-action2}:
\begin{cor}
  \label{cor:trivial}
  In the context above, for a group homomorphism $\alpha: G \to
  \mathrm{Aut}(\Pi)$ given by \autoref{eq:20} with associated isomorphism \autoref{eq:19}, 
the following categories are equivalent: 
\begin{itemize}
\item A $G$-fiber bundle with fiber $F$ and structure group $\Pi$;
\item A $G$-fiber bundle with fiber $F'$, structure group $\Pi$ and total group $\Pi \times G$;
\item A $G$-fiber bundle with fiber $F'$, structure group $\Pi$ and total group $\Pi
  \rtimes_{\alpha} G$. \qed
\end{itemize}
\end{cor}
\noindent
Similarly, a principal $(\Pi; \Pi \times G)$-bundle is literally the same thing as a principal $(\Pi;\Pi
\rtimes_{\alpha} G)$-bundle, but they have different canonical $G$-actions.
\begin{notn}
  \label{notn:changeOfG}
  For a principal $G$-$\Pi$-bundle, we call it a principal $(\Pi; \Pi \times G)$-bundle if we let
  $G$ act on the total space by $G \subgroup \Pi \times G$; we call it a principal $(\Pi; \Pi \rtimes
  G)$-bundle if we let $G$ act on the total space by $\Lambda_{\beta} \subgroup \Pi \times G$.
  And similarly for a $G$-fiber bundle with fiber $F$ and structure group $\Pi$.
\end{notn}

This trivial observation allows us to define and study the $V$-framing bundle of an equivariant
vector bundle.
  Let $V$ be an orthogonal $G$-representation given by
  $\rho: G \to O(n)$. In the remainder of this subsection, we write $O(V)$ for the group $O(n)$ with the data $G
  \to \mathrm{Aut}(O(n))$ given by $g(\ele) = \rho(g) \ele \rho(g)^{-1}$ for $g \in G$ and $\ele \in
  O(n)$, so it is clear what $O(V) \rtimes G$ means.
  This convention coincides with the conjugation $G$-action
  on $O(V)$ thought of as a mapping space in $\topG$.
  In this case, taking $F=\bR^n$ and pointing aloud the $G$-action on $F'$, \autoref{cor:trivial}
  reads: A $G$-$n$-vector bundle is a $G$-fiber bundle with fiber $\bR^n$, structure group
  $O(n)$ and total group $O(n) \times G$, as well as a $G$-fiber bundle with fiber $V$, structure group
  $O(n)$ and total group $O(V) \rtimes G$.
  \begin{defn}
    \label{defn:frv}
  Let $p:E \to B$ be a $G$-$n$-vector bundle.
  Let $\mathrm{Fr}_V(E)$ be the space of the admissible maps with the $G$-action
$  g(\psi) = g \psi \rho(g)^{-1}.$
\end{defn}
\noindent In other words, $\mathrm{Fr}_V(E)$ has the same underlying space as
$\mathrm{Fr}_{\bR^n}(E)$, but we think of admissible maps as mapping out of $V$ instead of
$\bR^n$.
\begin{prop} \label{prop:frv}
  $\mathrm{Fr}_V(E)$ is a principal $(O(n);O(V) \rtimes G)$-bundle and
  we have isomorphisms of $G$-vector bundles:
\begin{equation*}
E \cong (\mathrm{Fr}_V(E) \times V)/O(n).
\end{equation*}
\end{prop}
\begin{proof}
  This is a corollary of the structure theorem \autoref{thm:G-structure-2}.
  Namely, \autoref{cor:trivial} and the
  explanation afterwards have turned the vector bundle $p: E \to B$
  into a $G$-fiber bundle with fiber $V$, structure group
  $O(n)$ and total group $O(V) \rtimes G$. 
  By examination, $\mathrm{Fr}_V(E)$ in \autoref{defn:frv}
  agrees with the construction $\mathrm{Fr}_V(E)$ in \autoref{thm:G-structure-2}. 
\end{proof}

\subsection{Fixed point theorems}
\label{sec:fixed-point-theorems}
  Non-equivariantly, the long exact sequence of the homotopy groups of a fiber sequence is a useful
  tool to study the homotopy group of one term, knowing the other two. To do this equivariantly,
  we need to know what taking-fixed-points does to equivariant bundles.
  We focus on $\Gamma = \Pi \times G$ in this subsection.

  Let $\mathrm{Rep}(G,\Pi)$ be the set:
\begin{equation*}
\mathrm{Rep}(G,\Pi)  = \{ \text{group homomorphism }\rho: G \to \Pi\}/ \Pi
 \text{-conjugation}.
\end{equation*}

Any subgroup $H \subgroup G$ with a group homomorphism $\rho: H \to \Pi$ gives a subgroup $\Lambda_{\rho}$ of
$(\Pi \times G)$ via its graph. That is,
   $$\Lambda_{\rho} = \{(\rho(h) ,h)| h \in H\}.$$
   
For each $\rho: H \to \Pi$, denote  the centralizer
  of the image of $\rho$ in $\Pi$ by
\begin{equation*}
\mathrm{Z}_{\Pi}(\rho) = \{\ele \in \Pi| \ele \rho(h) = \rho(h) \ele  \text{
    for all } h \in H\}.
\end{equation*}

  \begin{prop}
    Let $\Pi$ be a compact Lie group and $H$ be a subgroup.
    Then $\mathrm{Z}_{\Pi}(H)$ is a closed subgroup of $\Pi$, thus also a compact Lie group.
  \end{prop}
  \begin{proof}
    Fix an element $h \in H$. Then the map $c_h:\Pi \to \Pi, \ \ele \mapsto \ele h \ele^{-1}$ is
    continuous. Since the singleton $\{h\} \in \Pi$ is closed, the set $c_h^{-1}(\{h\}) =  \{\ele \in \Pi|
    \ele h = h \ele\}$ is also closed. So $\mathrm{Z}_{\Pi}(H) = \bigcap_{h \in H} c_h^{-1}(\{h\})$ is
    closed.
  \end{proof}

  The following theorem is proven in \cite[Theorem 12]{LM86} for
  general~$\Gamma$. We spell it out for $\Gamma = \Pi \times G$.
  \begin{thm}
    \label{thm:LM}
    Let $G$ and $\Pi$ be compact Lie groups.
    Let ${p: E \to B}$ be a principal $G$-$\Pi$-bundle and $H\subgroup G$ be a subgroup. Assume that $E$ is
    completely regular.
\begin{enumerate}
\item \label{item:LM1} On the base,
   \begin{equation*}
     B^H = \coprod_{[\rho] \in \mathrm{Rep}(H,\Pi)} p(E^{\Lambda_{\rho}}).
\end{equation*}
\item \label{item:LM2} As sets, the preimages over each component of $B^{H}$ are
\begin{equation*}
p^{-1}(p(E^{\Lambda_{\rho}})) = \coprod_{\{\rho': \Pi \text{-conjugate to } \rho \}} E^{\Lambda_{\rho'}}.
\end{equation*}
As spaces,
\begin{equation*}
p^{-1}(p(E^{\Lambda_{\rho}})) \cong \Pi \times_{Z_{\Pi}(\rho)} E^{\Lambda_{\rho}}.
\end{equation*}
\item \label{item:LM3} For a fixed representative $\rho$ of $[\rho]$, we have a principal
$Z_{\Pi}(\rho)$-bundle:
\begin{equation*}
Z_{\Pi}(\rho) \to E^{\Lambda_{\rho}} \overset{p}{ \to } p(E^{\Lambda_{\rho}}).
\end{equation*}
\item \label{item:LM4} In particular, the following is a principal $\Pi$-bundle:
\begin{equation*}
  \Pi \to E^{H} \overset{p}{ \to } p(E^{H}).
\end{equation*}
\end{enumerate}
\end{thm}
\begin{proof}[Explanation]
   In words, part \autoref{item:LM1}
   says that the $H$-fixed points of $B$ are the images of the {$\Lambda$-fixed} points
   of $E$ for all subgroups $\Lambda \subgroup  \Pi \times G$ that are graphs of a
   homomorphism $H \to \Pi$.
   Furthermore,  $E^{\Lambda}$ and $E^{\Lambda'}$ share the same projection image when
   $\Lambda$ and $\Lambda'$ are $\Pi$-conjugate, or equivalently the corresponding
   representations $H \to \Pi$ are $\Pi$-conjugate.
   The assumption that $E$ is completely regular implies that if $\Lambda$ and $\Lambda'$ are
   not $\Pi$-conjugate, the images of $E^{\Lambda}$ and $E^{\Lambda'}$ are disjoint.

   Parts \autoref{item:LM2} and \autoref{item:LM3} imply that $E$ restricted on each component of
   $B^H$ has a reduction of the structure group from $\Pi$ to $Z_{\Pi}(\rho)$.
   In the proof of \autoref{thm:BGO(n)2}\autoref{item:BGO(n)fiber}, we will describe in an example how to find
   the representations $\rho$ when $H = G$. The idea is that the fiber over an $H$-fixed base
   has an $H$-action, and $\rho$ tells what this action is in terms of the native $\Pi$-action
   as a principal bundle.
   Note that the representation $\rho$ is dependent on the choice of a base point $z$ in the fiber;
   a different choice gives a conjugate representation.
   From the description of the action, 
   a point in the same fiber, written uniquely as $z \ele$ for some $\ele \in \Pi$,
   is $\Lambda_{\rho}$-fixed if and only if $\rho(h) \ele \rho(h)^{-1}=\ele$ for all $h \in H$.
   This justifies the first statement of part \autoref{item:LM2} as well as part \autoref{item:LM3}.

   For the second statement of part \autoref{item:LM2}, which is not in the reference, we use the
   map:
\begin{equation*}
\Pi \times_{Z_{\Pi}(\rho)} E^{\Lambda_{\rho}} \to E, \ (\ele,x) \mapsto x\ele^{-1}.
\end{equation*}
Here, $Z_{\Pi}(\rho)$ is a subgroup of $\Pi$ and acts on the right of $\Pi$ by multiplication; the
left $\Pi$-action on $E$ restricts to a left $Z_{\Pi}(\rho)$-action on $E^{\Lambda_{\rho}}$ .
It is a homeomorphism to its image, which is exactly $p^{-1}(p(E^{\Lambda_{\rho}}))$:

We have $\Lambda_{e} = H$ for the trivial representation $e: H \to \Pi$. 
Part \autoref{item:LM4} follows from taking $\rho = e$ in part \autoref{item:LM3}.
\end{proof}

\begin{rem}
  \label{rem:obeservation-rep}
  From \autoref{thm:LM}, for a principal $G$-$\Pi$-bundle $p:E \to B$ and a subgroup $H \subgroup G$,
  each component $B_0$ of $B^H$ has an associated representation class
  $[\rho] \in \mathrm{Rep}(H,\Pi)$. It is characterized by the fact that  
  for any representation $\rho': H \to \Pi$,
  \begin{equation*}
    \big(p^{-1}(B_0)\big)^{\Lambda_{\rho'}} \not= \varnothing \text{ if and only if } [\rho'] = [\rho].
\end{equation*}
The restricted principal $\Pi$-bundle
$p^{-1}(B_0) \to B_0$ has a reduction of the structure group from $\Pi$ to $Z_{\Pi}(\rho)$.
\end{rem}

  Non-equivariantly, a map between two principal $G$-bundles that is an underlying
  equivalence on the total spaces will give an equivalence on the base spaces,
  as can be shown by the long exact sequence of homotopy groups.
  Equivariantly, we also want this tool of knowing when a map of two principal $G$-$\Pi$-bundles
  gives a $G$-equivalence on the base spaces.
   \begin{thm} \label{lem:comparefib} 
    Let $G$,$\Pi$ be compact Lie groups and $i: \Pi \to \Pi' $ be an inclusion
    of topological groups. Let $p:E \to B$ be a principal $G$-$\Pi$-bundle, $p':  E' \to B'$ be a principal $G$-$\Pi'$-bundle and assume that $B,B'$ have the homotopy
    type of $G$-CW complexes.
    Then $E'$ has a $(\Pi \times G)$-action by $i$.
    
    Suppose that there is a $(\Pi \times G)$-map $\bar{f}: E \to E'$ over a $G$-map $f: B \to B'$,
    as in the following commutative diagram:
    \begin{equation*}
      \begin{tikzcd}
        \Pi \ar[r,"i"] \ar[d] & \Pi' \ar[d]\\
        E \ar[r,"\bar{f}"]  \ar[d,"p"] & E'  \ar[d,"p'"] \\
        B \ar[r,"f"] & B' 
      \end{tikzcd}
    \end{equation*}
Assume
\begin{enumerate}[(1)]
\item \label{item:condition1} The map $i$ includes $\Pi$ as a deformation retract of $\Pi'$ in groups, that is, there exists a group
  homomorphism $j: \Pi'  \to \Pi$ such that $j\circ i = \mathrm{id}$ and $i
  \circ j \simeq \mathrm{id} \text{ rel }i(\Pi)$ in topological groups;
\item \label{item:condition2} On the total spaces, the map $\bar{f}$ is an
  $\mathscr{F}(\Pi)$-equivalence, that is, a $\Lambda$-equivalence for any subgroup
  $\Lambda \subset \Pi \times G$ such that $\Lambda \cap \Pi = e$.
\end{enumerate}
    Then, on the base spaces, $f: B \to B'$ is a $G$-equivalence.
  \end{thm}

  \begin{proof}
    To simplify notation in this proof, we use the same letters to denote the restrictions of
    the corresponding maps to a subspace.    
    By the equivariant Whitehead theorem, it suffices to show that:
\begin{equation*}
   \text{For any subgroup }H\subgroup G, \text{ the map } f: \ B^H \to (B')^H \text{  is an equivalence.}
\end{equation*}

We make the following two claims comparing $\Pi$ and $\Pi'$:
\begin{enumerate}[(a)]
\item \label{rmk:group(2)}For any group $H$, the induced map $i_{*} : \mathrm{Rep}(H,\Pi)
     \to \mathrm{Rep}(H, \Pi') $ is a bijection.
\item \label{rmk:group(1)}For any subgroup $K$ of $\Pi$, the inclusion $i : Z_\Pi K \to
  Z_{\Pi'} i(K)$ is a homotopy equivalence;
\end{enumerate}
These two claims follow from the assumption \autoref{item:condition1}.
For \autoref{rmk:group(2)}, we take the functor $F=\mathrm{Rep}(H,-)$ from the category of groups to
sets. It has equivalent images on $\Pi$ and $\Pi'$, and we skip the details.
For \autoref{rmk:group(1)}, we take the functor $F=Z_{(-)}K$ from the category of groups containing
$K$ as a subgroup. It also has equivalent images on $\Pi$ and $\Pi'$, and the details come later in
\autoref{lem:alemma}.
   
    By \autoref{thm:LM}~\autoref{item:LM1} and \autoref{rmk:group(2)},
    it suffices to show that:
\begin{equation*}
   \text{For any $H$ and }\rho \in \mathrm{Rep}(H,\Pi),\text{ the map }f: p(E^{\Lambda_{\rho}}) 
   \to  p'((E')^{\Lambda_{\rho}}) \text{  is an equivalence.}
\end{equation*}
By \autoref{thm:LM}~\autoref{item:LM3},
taking the $\Lambda_\rho$-fixed points of $E$ and $E'$ yields a map between principal bundles:
    \begin{equation*}
      \begin{tikzcd}
        Z_{\Pi}(\rho) \ar[r,"i"] \ar[d] & Z_{\Pi'}(\rho) \ar[d] \\
        E^{\Lambda_{\rho}} \ar[r,"\bar{f}"] \ar[d,"p"] & (E')^{\Lambda_{\rho}} \ar[d,"p'"] \\
        p(E^{\Lambda_{\rho}}) \ar[r,"f"] & p'((E')^{\Lambda_{\rho}})
      \end{tikzcd}
    \end{equation*}
    By the claim \autoref{rmk:group(1)} and the assumption \autoref{item:condition2},
    both $i$ and $\bar{f}$ are equivalences. The long exact
    sequence of homotopy groups shows that $f$ is an equivalence.
  \end{proof}

  \begin{rem}
    \label{rmk:group}
    In \autoref{lem:comparefib}, the assumption \autoref{item:condition1}
    is satisfied in our applications with $\Pi' = \Pi$ or $\Pi' = \Pi^{I}$.    
    The assumption \autoref{item:condition2} is
    satisfied when $\bar{f}$ is a $(\Pi \times G)$-equivalence, but is weaker.
    The weaker version is needed in our applications.
  \end{rem}

  From the proof, we also have a version of \autoref{lem:comparefib} relaxing the
  assumption~\autoref{item:condition2}.
  \begin{cor}
    \label{cor:obeservation-comparefib}
    Suppose we have $(i, \bar{f}, f)$ in the context of \autoref{lem:comparefib}, except that
    instead of the assumption \autoref{item:condition2}, $\bar{f}: E \to E'$
    is only known to be a $\Lambda_{\rho}$-equivalence
    for a fixed representation $\rho: H \to \Pi$. Then on the base spaces,
    $f: p(E^{\Lambda_{\rho}}) \to p((E')^{\Lambda_{\rho}})$ is an equivalence.
  \end{cor}
  Note that $p(E^{\Lambda_{\rho}})$ is the space of components of $B^H$ that are associated to $\rho$ as
  described in \autoref{rem:obeservation-rep}. In particular, if $(B')^H$ is connected for all
  subgroups $H \subgroup G$, then $(B')^H$ has only one associated representation $\rho_H$. Moreover,
  $\rho_H$ has to be the restriction of $\rho_G$. We have:
    \begin{cor}
      Let $B'$ be a $G$-connected space as explained above and $\rho_G$ be the associated
      representation.
      Suppose we have $(i, \bar{f}, f)$ in the context of \autoref{cor:obeservation-comparefib},
     such that $\bar{f}$ is a $\Lambda_{\rho_G}$-equivalence.
     Then on the base spaces,  $f: B \to B'$ is a $G$-equivalence.
    \end{cor}
    \begin{proof}
      Since the map $f: B^H \to (B')^H$ preserves the associated representation, we know that $B^H$
      only has one associated representation $\rho_H$ as well. The claim then follows by
      applying \autoref{cor:obeservation-comparefib} to $\rho = \rho_H$ for all $H$.
    \end{proof}
 
   The following is a lemma for \autoref{lem:comparefib}:
 \begin{lem}
      \label{lem:alemma}
   Assume $i:\Pi \to \Pi'$ is an inclusion of topological groups with a deformation retract $j:\Pi'
   \to \Pi$, that is, they satisfy condition \autoref{item:condition1} in \autoref{lem:comparefib}.
   Then for any  subgroup $K$ of $\Pi$, the
   inclusion $i : Z_\Pi K \to  Z_{\Pi'} i(K)$ is a homotopy equivalence.
 \end{lem}
 \begin{proof}
   We first check that in general, given any group homomorphism $f: G \to G'$ and subgroup $K
   \subgroup G$, the map $f$ restricts to a map 
   $f_0: Z_GK \to Z_{G'}(f(K))$ on subspaces. This is because $x k = k x$ for all $k \in K$ implies
   $f(x)f(k)=f(k)f(x)$ for all $f(k) \in f(K)$.
   So, we have 
\begin{equation*}
i_0:  Z_\Pi K \to  Z_{\Pi'} (i(K)) \text{ and } j_0: Z_{\Pi'} (i(K)) \to  Z_\Pi (ji(K)) = Z_\Pi K.
\end{equation*}
The map $j_0$ gives deformation retract data of the inclusion $i_0$.
It is obvious that $j_0i_0 = \mathrm{id}$. It remains to show $i_0j_0 \simeq \mathrm{id}$.
The image of $i_0$ is the subspace $Z_{i(\Pi)}(i(K)) \subset Z_{\Pi'}(i(K))$. The homotopy $ij \simeq
\mathrm{id}$ rel $i(\Pi)$ restricts to a homotopy $i_0j_0 \simeq \mathrm{id}$ rel $Z_{i(\Pi)}(i(K))$.
 \end{proof}

\section{Classifying spaces}
\label{sec:classify}

\subsection{$V$-trivial bundles}
\label{sec:v-trivial-bundles}
An equivariant bundle $E \to B$ is $V$-trivial for some $n$-dimensional $G$-representation
$V$ if there is a $G$-vector bundle isomorphism $E \cong B \times V$.
Such an isomorphism is a $V$-framing of the bundle. This is analogous to the
case of non-equivariant vector bundles, except that equivariance adds in the complexity
of a representation $V$ that's part of the data.

However, the representation $V$ in the equivariant trivialization of a fixed vector bundle may not be unique. 
We give a lemma to recognize when two trivial bundles are isomorphic, then a counterexample.

Let $\mathrm{Iso}(V,W)$ be the space of linear isomorphisms $V \to W$ with the conjugation
$G$-action for $G$-representations $V$ and $W$.
  \begin{lem}
\label{lem:trivialbundle}
    For a $G$-space $B$, there exists a $G$-vector bundle
    isomorphism $B \times V \cong B \times W$ if and only if there exists a $G$-map $f: B \to
    \mathrm{Iso}(V,W)$.
  \end{lem}
  \begin{proof}
    Let $F: B \times V \to B \times W$ be a vector bundle map. For $b \in
    B$, let $F_b: V \to W$ be such that $F_b(v) = F(b,v)$. Then $F$ is a $G$-vector bundle
    isomorphism if and only if
    \begin{enumerate}
    \item $F$ is fiberwise isomorphism: $F_b \in \mathrm{Iso}(V,W)$;
    \item $F$ is a $G$-map: $gF(b,v) = F(gb,gv)$, or equivalently,
      $F_{gb} = gF_bg^{-1}$, for all $g \in G$.
    \end{enumerate}
    Taking $f(b) = F_b$, it follows that $F$ is an isomorphism if and only if
    $f$ is a $G$-map. 
  \end{proof}
  \begin{cor}
\label{cor:trivialbundle}
    If $B$ has a $G$-fixed point, then $B \times V \cong B \times W$ only when $V \cong W$.
  \end{cor}
  \begin{proof}
    The equivariant map $f: B \to \mathrm{Iso}(V,W)$ induces $f^G: B^G \to
    \mathrm{Iso}_G(V,W)$. The source being nonempty implies that the target is nonempty.
  \end{proof}

  \begin{rem}
      More generally, for any two $n$-dimensional $G$-vector bundles $E,E'$ over
  $B$, one can form the non-equivariant bundle $\mathcal{H}om_{B}(E,E')$ which consists of all bundle
  maps $E \to E'$ over $B$  (not necessarily fiberwise isomorphisms).
  It has a $G$-action by conjugation and is indeed an $n^2$-dimensional $G$-vector bundle
  over $B$. Let $\mathcal{I}so_{B}(E,E')$ be the subspace consisting of only
  fiberwise isomorphisms. It is a $GL_n$-bundle over $B$. Then tautologically $E \cong E'$ if there is a
  $G$-invariant section of $\mathcal{I}so_{B}(E,E')$.
\end{rem}

  \begin{exmp}[Counterexample]
    Let $G=C_2$, $\sigma$ be the sign representation. The unit sphere, $S(2\sigma)$, is $S^1$ with the
  180 degree rotation action. As $C_2$-vector bundles,  
\begin{equation*}
S(2\sigma) \times \mathbb{R}^2 \cong S(2\sigma) \times 2\sigma.
\end{equation*}

  \begin{proof}
    By \autoref{lem:trivialbundle}, it suffices to construct a $C_2$-map ${S(2\sigma) \to
  \mathrm{Iso}(\mathbb{R}^2, 2\sigma) \cong GL_2}$, where the nontrivial element
  of $C_2$ acts on $GL_2$ by multiplying by $-\mathrm{Id}$.
  We give $S(2\sigma)$ a $G$-CW decomposition of a 0-cell
  $C_2/e$ and a 1-cell $C_2/e \times \mathrm{D}^1$ and construct the map by skeleton. It is obvious that any
  equivariant map on the 0-skeleton extends to the 1-skeleton if and only if the two images lie in the same
  path component of $GL_2$, which is true in this case as $-\mathrm{Id}$ and
  $\mathrm{Id}$ lie in the same path component.
  \end{proof}
  \end{exmp}

  The following counterexample is suggested by Gus Longerman.
  \begin{exmp}(Counterexample)
    Take $G$ to be any compact Lie group and $V$ and $W$ to be any two representation
    of $G$ that are of the same dimension. Then $G \times V  \cong G \times W$,
    because $\mathrm{Map}_G(G, \mathrm{Iso}(V,W)) \cong \mathrm{Map}(\mathrm{pt},
    \mathrm{Iso}(V,W)) \neq \varnothing$.
    Indeed, the isomorphism can be constructed explicitly by $F(g,x) = (g, \rho_W(g)
    \rho_V(g)^{-1}x)$, where $\rho_V , \rho_W: G \to O(n)$ are matrix
    representations of $V,W$.
  \end{exmp}

  \subsection{Universal equivariant bundles}
  \label{sec:univ-bundle}
\begin{defn}\label{thm:universalmap}
  A principal $(\Pi;\Gamma)$-bundle $E \to B$ is called universal if
  for any paracompact $G$-space $X$, there is a bijection of sets between \{equivalence classes of principal $(\Pi;\Gamma)$-bundles over $X$\}
and \{$G$-homotopy classes of $G$-maps $X \to B\}$.
\end{defn}
The correspondence in one direction is by pulling back the universal bundle along $G$-maps $X \to
B$. The universal principal $(\Pi;\Gamma)$-bundle exists and is unique up to
homotopy. Moreover, it can be recognized by the following property:

\begin{thm}(\cite[Theorem 9]{LM86})
\label{thm:universalbundle}
  A principal $(\Pi; \Gamma)$-bundle $p:E \to B$ is universal if and only if
\begin{equation*}
E^{\Lambda} \simeq *, \text{ for all subgroups } \Lambda \subset \Gamma 
\text{ such that } \Lambda \cap \Pi
= {e}.
\end{equation*}
\end{thm}

The universal principal $(\Pi; \Gamma)$-bundle can be constructed
    in various ways. We denote it as $E(\Pi;\Gamma) \to B(\Pi;\Gamma)$.
    One construction generalizes Milnor's infinite join construction. See
    \cite[Section 3]{TD} for the case $\Gamma = \Pi \rtimes G$ or
    \cite[Section 2]{L82} for the case $\Gamma = \Pi \times G$.
    We describe another abstract construction using
    \autoref{thm:universalbundle} as outlined in the introduction. Suppose we have a universal $\Gamma$-space  $E \mathscr{F}$ for the family
    $\mathscr{F} = \{\Lambda \subset \Gamma | \Lambda \cap \Pi = e \}$. Then the restricted $\Pi$-action on $E \mathscr{F}$ is free, as any non-trivial
    element $A \in \Pi$ will generate a subgroup $\langle A \rangle \subset \Gamma$ not in
    $\mathscr{F}$, so $\{x \in E \mathscr{F}| Ax = x\} \cong E
    \mathscr{F}^{\langle A \rangle} = \varnothing$. So the quotient map $E \mathscr{F}
    \to E \mathscr{F} / \Pi $  is a principal $(\Pi;\Gamma)$-bundle, and it is
    equivalent to the universal one by \autoref{thm:universalbundle}. Thus it suffices to construct $E
    \mathscr{F}$. This can be done via Elmendorf's construction
    \cite{Elmendorf}: The $E \mathscr{F}$ is constructed as a $\Gamma$-CW complex (See
    \autoref{sec:g-cw-complexes} for $G$-CW complexes).
    First, we take $$E_0 =
    \sqcup_{\Lambda \in \mathscr{F}} \Gamma/\Lambda.$$
    It has the property that $E_0^{\Lambda} \neq \varnothing$ exactly when $\Lambda \in
\mathscr{F}$. We make the desired fixed-point-spaces contractible by adding in
higher cells. To start, we add 1-cells 
\begin{equation*}
 \sqcup_{\Lambda_0, \Lambda_1 \in \mathscr{F}} \mathrm{Map}_\Gamma(\Gamma/\Lambda_1, \Gamma / \Lambda_0) \times
\Gamma/\Lambda_1 \times D^1
\end{equation*}
to obtain $E_1$. Here,
$\mathrm{Map}_\Gamma(\Gamma/\Lambda_1, \Gamma / \Lambda_0)$ is the space of
$\Gamma$-equivariant maps topologized as homeomorphic to the space $(\Gamma /
\Lambda_0)^{\Lambda_1}$. The gluing map on 
\begin{equation*}
(f,x, t) \in  \mathrm{Map}_\Gamma(\Gamma/\Lambda_1, \Gamma / \Lambda_0) \times
\Gamma/\Lambda_1 \times S^0
\end{equation*}
sends $(f,x,-1)$ to $f(x) \in \Gamma/\Lambda_0
\subset E_0$ and sends
$(f,x,1)$ to $x \in \Gamma/\Lambda_1 \subset E_0$. Now we have $E_1^{\Lambda}$
is connected for $\Lambda \in \mathscr{F}$. The higher cells are selected and glued following the same
idea and we skip the details. Technically speaking, when $\Gamma$ is not a discrete group, $ \mathrm{Map}_\Gamma(\Gamma/\Lambda_1, \Gamma / \Lambda_0) \times
\Gamma/\Lambda_1 \times D^1$ may not be a 1-cell as defined in
\autoref{sec:equiBundle}. However, since $\mathrm{Map}_\Gamma(\Gamma/\Lambda_1, \Gamma
/ \Lambda_0)$ has the homotopy type of a CW-complex, the constructed $E
\mathscr{F}$ always has
the homotopy type of a $\Gamma$-CW complex.

\begin{rem}
  \label{rem:Lambeda}
When $\Gamma = \Pi \times G$, such a subgroup $\Lambda$ comes in the form of 
\begin{equation*}
\{(\rho(h),h) | h \in H \}, \text{ for } H \subgroup  G \text{ and } \rho : H \to \Pi \text{ is a group homomorphism}.
\end{equation*}
This group was denoted by $\Lambda_{\rho}$ in \autoref{thm:LM}.

When $\Gamma = \Pi \rtimes_{\alpha} G$, such a subgroup $\Lambda$ comes in the form of 
\begin{equation*}
\{(\rho(h),h) | h \in H \}, \text{ for }H \subgroup  G \text{ and } \rho : H \to \Pi
\text{ such that } \rho(h_1h_2) = \rho(h_1) \cdot \alpha_{h_1}(\rho(h_2)).
\end{equation*}
\end{rem}

We mostly specialize to the case $\Gamma = O(n) \times G$, when a principal $(O(n);\Gamma)$ is also a
principal $G$-$O(n)$-bundle. We also denote the universal principal
$G$-$O(n)$-bundle by $E_GO(n) \to B_GO(n)$ and denote the universal $G$-$n$-vector bundle by $\universal \to B_GO(n)$
where $$\universal = E_GO(n) \times_{O(n)} \bR^n.$$

\medskip

As an immediate corollary of Theorems~\ref{thm:LM} and \ref{thm:universalbundle}, one gets the
  $G$-homotopy type of the universal base.
Recall that
  \begin{align*}
 \mathrm{Rep}(G,O(n)) & = \{\rho: G \to O(n) \text{ group homomorphism }\}/{O(n)\text{-conjugation}}; \\
   & \cong \{V: n \text{-dimensional orthogonal representation of } G\}/\text{isomorphism}
  \end{align*}
  and  $\mathrm{Z}_{O(n)}(\rho) = \{a \in O(n)| a \rho(g) = \rho(g) a,  \text{
    for all } g \in G\}$ is the centralizer
  of the image of $\rho$ in $O(n)$.

\begin{thm}(\cite[Theorem 2.17]{L82})
\label{thm:BGO(n)}
\begin{align*}
(B_GO(n))^{G} & \simeq \coprod_{[\rho] \in
                \mathrm{Rep}(G,O(n))}B\mathrm{Z}_{O(n)}(\rho); \\
    & \simeq \coprod_{[V] \in
                \mathrm{Rep}(G,O(n))}B (O(V)^G).
\end{align*}
 \end{thm}

\begin{exmp}
  Take $H=G=C_2$ and $\Pi=O(2)$. Then
  $$\mathrm{Rep}(C_2,O(2)) = \{\mathrm{id}, \text{ rotation}, \text{
    reflection}\},$$
  where $\mathrm{rotation}$ is the rotation by 180 degrees and
  $\mathrm{reflection}$ is the equivalence class of all reflections.
  For $\rho=\mathrm{id}$ or $\rho=\text{rotation}$, $Z_{\Pi}(\rho)=O(2)$. For
  $\rho=\text{reflection}$, $Z_{\Pi}(\rho) \cong \bZ/2 \times \bZ/2$. So 
\begin{equation*}
(B_{C_2}O(n))^{C_2} \simeq BO(2) \sqcup BO(2) \sqcup B(\bZ/2 \times \bZ/2).
\end{equation*}
\end{exmp}

\subsection{The fiber of the universal bundle and the loop space of the
  classifying space over a point}
\label{sec:fiber-univ-bundle}
  From \autoref{thm:BGO(n)}, one can make explicit the classifying maps of
  $V$-trivial bundles as follows.
  
   A $G$-map $\theta: \mathrm{pt} \to B_GO(n)$ lands in one of the $G$-fixed
   components of $B_GO(n)$. Suppose that this component is indexed by $[V]$. 
  \begin{prop}
\label{rmk:classifyingV}
 The pullback of the universal bundle
  is $\theta^{*}\universal \cong V$ over $\mathrm{pt}$.
\end{prop}
\begin{proof}
  It follows from part \autoref{item:BGO(n)fiber} of the following
  \autoref{thm:BGO(n)2} that $$\theta^{*} \universal \cong O(\bR^n,V) \times_{O(n)} \bR^n \cong V.$$
  In fact, the $n$-planes in a complete $G$-universe with the tautology $n$-plane bundle is a model
  for $B_GO(n)$ and $\universal$; $\theta(\mathrm{pt})$ is just a $G$-representation isomorphic to $V$.
\end{proof}

  \begin{thm}
\label{thm:BGO(n)2}
      Take a $G$-fixed base point $b \in B_GO(n)$ in the component indexed by $[V]$. Let $p:E_GO(n) \to
      B_GO(n)$ be the universal principal $G$-$O(n)$-bundle. Then
\begin{enumerate}
\item\label{item:BGO(n)fiber} The fiber over $b$, $p^{-1}(b)$, is homeomorphic to
  $O(\mathbb{R}^n, V)$ as an $(O(n) \times G)$-space. Here, $(
  O(n) \times G)$ acts on $O(\mathbb{R}^n, V)$ by $G$ acting on $V$ and $O(n)$
  acting on $\mathbb{R}^n$.
\item\label{item:BGO(n)loop} The loop space of $B_GO(n)$ at the base point $b$,
  $\Omega_b B_GO(n)$, is $G$-homotopy equivalent to $O(V)$, the isometric self maps of $V$ with $G$
  acting by conjugation.
\end{enumerate}
  \end{thm}
  \begin{proof}

    \autoref{item:BGO(n)fiber} This is due to Lashof and we explain how to find the representation
    $V$ here. Choose and fix a base point $z \in p^{-1}(b)$. We
construct a group homomorphism $\rho_{z}: G \to O(n)$ as follows. For any
$g \in G$, there exists a unique element,  $\rho_{z}(g) \in O(n)$ such
    that $g z = z \rho_{z}(g)$. It is easy to check that $g \mapsto \rho_{z}(g)$
    gives a group homomorphism.  Suppose $z'$ is another base point in
    $p^{-1}(b)$, and $z'=z\ele$ for some unique $\ele \in O(n)$. Then
\begin{equation*}
    gz' = gz\ele = z \rho_z(g) \ele = z'(\ele^{-1}\rho_z(g)\ele).
\end{equation*}
    So $\rho_{z'} = \ele^{-1}\rho_z\ele$ is $O(n)$-conjugate to $\rho_z$. 
    The different $\rho_{z}$'s are the matrix representations of some vector space representation $V$.
    From the proof of Theorem 2.17 of \cite{L82}, this is exactly the index
    $V$. Without loss of generality, we take $V$ to be given by $\rho_z$ as matrix representation.

    The following map gives a non-equivariant homeomorphism:
    \begin{equation*}
      \begin{array}{cccc}
       O(\mathbb{R}^n,V) & \cong O(n)& \overset{\cong}{ \to } & p^{-1}(b), \\
           & \ele & \mapsto & z\ele.
      \end{array}
    \end{equation*}
   It suffices to check it is an equivariant homeomorphism with the described action. Let $(\elm,g) \in
   O(n) \times G$. Then
    \begin{equation*}
      z((\elm,g) \circ \ele) = z (\rho_{z}(g)\ele \elm^{-1}) =(z \rho_z(g))(\ele\elm^{-1})
      =(gz)(\ele\elm^{-1}) = (\elm,g) \circ z\ele.
    \end{equation*}

\autoref{item:BGO(n)loop} The idea is to compare the path space fibration with the universal bundle.
Equivariantly, the base point should be $G$-fixed. Since the space involved is not $G$-connected,
base points from different components might give inequivalent loop spaces.
We use subscripts in path spaces and loop spaces to indicate the base point. For example,
$$P_{b}B_GO(n) = \{ \alpha \in \mathrm{Map}([0,1],B_GO(n)) | \alpha(0)=b \}.$$

Fix $z \in p^{-1}(b)$ and $\rho = \rho_{z}: G \to O(n)$ as above.
Take $z$ to be the base point of $E_GO(n)$. It is a $\Lambda$-fixed point, where
$$\Lambda = \{(\rho(g), g)| g \in G\} \subset O(n) \times G.$$
We prove that  $E_GO(n)$ is $\Lambda$-contractible.
In fact, let $\Lambda'$ be any subgroup of $\Lambda$. Then $\Lambda' \cap O(n) = e$, so by
\autoref{thm:universalbundle}, $(E_{G}O(n))^{\Lambda'}$ is contractible.

So, the contraction map gives a based $\Lambda$-equivariant homotopy:
\begin{equation*}
E_GO(n) \wedge I \to E_GO(n).
\end{equation*}
Here, $I=[0,1]$ is based at 0 and has the trivial $\Lambda$-action.
(The map sends $x \sm 0$ and $z \sm t$ to $z$ for all $x \in E_GO(n)$ and $t \in I$.)
We take the adjoint of this homotopy to get $ E_GO(n) \to P_{z}E_GO(n)$,
and then compose with $ P_{z}E_GO(n) \to P_{b}B_GO(n)$ induced
by  $p: E_GO(n) \to B_GO(n)$. The composite is 
\begin{equation*}
f: E_GO(n) \to P_{z}E_GO(n) \to P_{b}B_GO(n).
\end{equation*}
It sends a point $x \in E_GO(n)$ to a path in $B_GO(n)$ that starts at $b$ and ends at $p(x)$.
This yields a commutative diagram:
\begin{equation}
  \label{eq:1}
\begin{tikzcd}
E_GO(n) \ar[r,"f"] \ar[d,"p"'] & P_{b}B_GO(n) \ar[d,"p_1"] \\
B_GO(n) \ar[r,equal] & B_GO(n)
\end{tikzcd}
\end{equation}
Moreover, this diagram is $G$-equivariant, where
the $G$-action on $P_bB_GO(n)$ is by pointwise action on the path.
It is worth noting that the $G$-action we take on $E_GO(n)$ is not the original one, but via the
identification $q: \Lambda \cong G$. In other words, $g \in G$ acts by what $(\rho(g),g)$ acts.
The two vertical maps are non-equivariant fibrations and
$f$ maps the fiber of $p$ over $b \in B_GO(n)$, denoted $F_1$, to the fiber of $p_1$ over $b$,
denoted $F_2$.

We first identify the fibers $F_1$ and $F_2$.
From part \autoref{item:BGO(n)fiber}, $F_1 \cong O(\bR^n, V)$ as $(O(n) \times G)$-spaces.
So $F_1 \cong O(V)$ as $G$-spaces.
It is clear that $F_2 \cong \Omega_bB_GO(n)$ as $G$-spaces.

We claim that $f$ restricts to a $G$-equivalence $F_1 \to F_2$.
The strategy is to show that it induces an isomorphism on homotopy groups of $H$-fixed points
for all $H\subgroup G$, using the long exact sequences of homotopy groups of fiber sequences.
Without dealing with general $G$-fibrations,
it suffices to work out the following:
\begin{itemize}
\item Denote by $\Lambda' = q^{-1}(H)$, the subgroup of
  $\Lambda$ that is isomorphic to $H$. The commutative diagram \autoref{eq:1} restricts to the following
  commutative diagram:
\begin{equation*}
  \begin{tikzcd}
    (F_1)^H \ar[d] \ar[r] & (F_2)^{H}  \ar[d] \\
    (E_GO(n))^{\Lambda'} \ar[d,"p"'] \ar[r,"f^H"] & (P_bB_GO(n))^{H} \ar[d, "p_1"]\\
    p((E_GO(n))^{\Lambda'}) \ar[r] & p_1((P_bB_GO(n))^{H})
  \end{tikzcd}
\end{equation*}
\item On the total space level, $f^H$ induces isomorphism on homotopy groups. This is true because
  $E_GO(n)$ is $\Lambda$-contractible and $P_bB_GO(n)$ is $G$-contractible.
\item The base spaces are equal. In fact, it is easy to see that they are both the component of
  $(B_GO(n))^H$ indexed by $[V]$ from Theorems \ref{thm:LM} and  \ref{thm:BGO(n)}.
\item The two vertical lines are fiber sequences. For the first, we use
  \autoref{thm:LM}~\autoref{item:LM3} with $(F_1)^H = (O(V))^H = Z_{\Pi}(\rho|_H)$;
  for the second, it is merely the path space fibration $\Omega_b X \to P_b X \to X$,
  where $X$ denotes the component of $(B_GO(n))^H$ containing $b$. \qedhere
\end{itemize}
\end{proof}

\begin{rem}
  \label{rem:general-Pi} The proof of \autoref{thm:BGO(n)2} works for general $\Pi$ placing
  $O(n)$.  Take a $G$-fixed base point $b \in B_G\Pi$ in the component indexed by
  $[\rho: G \to \Pi]$. Let $\Pi_{\mathrm{ad}}$ be the space $\Pi$ with the adjoint $\Pi$-action and
  consider it as a $G$-space via $\rho$. Then there is a $G$-homotopy equivalence
  $\Omega_b B_G\Pi \simeq \Pi_{\mathrm{ad}}$.
\end{rem}

\subsection{The gauge group of an equivariant principal bundle}
\label{sec:gauge-group}
Let $EO(n) \to BO(n)$ be the universal principal $O(n)$-bundle and ${p: P \to B}$ be any principal
$O(n)$-bundle. The gauge group of $P$, $\mathrm{Aut}_B(P)$, is the space of bundle
automorphisms of $P$ that are identity on the base space $B$ (\cite[Chap 7, Definition 1.1]{Husemoller}).
It turns out that the space of principal bundle maps, $\mathrm{Hom}(P, EO(n))$, is also universal:
The map 
\begin{equation}
  \label{eq:21}
\mathrm{Hom}(P, EO(n)) \to \mathrm{Map}_{p}(B, BO(n))
\end{equation}
that restricts a bundle map to its
base spaces is known to be the universal principal $\mathrm{Aut}_B(P)$-bundle.
Here, $\mathrm{Map}_p(B, BO(n))$ denotes the component of the
classifying map of $p$ in $\mathrm{Map}(B, BO(n))$.
A proof of this result can be found in \cite[Chap 7, Corollary 3.5]{Husemoller}.
In this subsection, we show the equivariant generalization of this result
(\autoref{thm:MapToUniversal}).

Let  $E_GO(n) \to B_{G}O(n) $ be the
universal principal $G$-$O(n)$-bundle and $p : P \to B$ be any principal $G$-$O(n)$-bundle.
The restricting-to-the-base map
\begin{equation}
  \label{eq:universal}
\pi: \mathrm{Hom}(P, E_GO(n)) \to \mathrm{Map}_p(B,B_GO(n))
\end{equation}
is a $G$-map lifting \autoref{eq:21}. Here, $\mathrm{Map}_p(B,B_GO(n))$ is the (non-equivariant)
component of the classifying map of $p$ in $\mathrm{Map}(B,B_GO(n))$; $G$ acts by conjugation on
both sides of \autoref{eq:universal}.
Let $\Gamma = \mathrm{Aut}_{B}P \rtimes G$, where $G$ acts on $\mathrm{Aut}_BP$ by conjugation.
Then the map $\pi$ in \autoref{eq:universal} is a universal principal
$(\mathrm{Aut}_B(P);\Gamma)$-bundle.
Note that this is an equivariant principal bundle not in the sense of \autoref{defn:Gprincipal}, but
of \autoref{defn:GGammaprincipal} - the total group is a non-trivial extension of
$\mathrm{Aut}_B(P)$ by $G$.

\begin{thm}
\label{thm:MapToUniversal}
In the context above,
the map 
\begin{equation*}
\pi: \mathrm{Hom}(P, E_GO(n)) \to \mathrm{Map}_p(B,B_GO(n))
\end{equation*}
 is a principal  $(\mathrm{Aut}_{B}P; \Gamma)$-bundle, where $\Gamma =
 \mathrm{Aut}_B(P) \rtimes G$. 
 Moreover, we have $$ \mathrm{Hom}(P, E_GO(n)) \simeq E \mathscr{F}$$ for the family $\mathscr{F} = \{\Lambda \subgroup  \Gamma \text{ such
   that } \Lambda \cap \mathrm{Aut}_BP = e\}$.
\end{thm}}
\begin{proof}
  As stated above, it is known non-equivariantly that $\pi$ is a universal principal
  $\mathrm{Aut}_BP$-bundle. One can use the conjugation $G$-action to get a principal
  $(\mathrm{Aut}_{B}P; \Gamma)$-bundle structure on $\pi$.
  However, later in this proof we want a $\Gamma$-action on the bundle $P$,
  so at the risk of elaborating the obvious, we describe the $\Gamma$-action on
  $\mathrm{Hom}(P,E_GO(n))$ by putting a $\Gamma$-action on both $P$ and $E_GO(n)$. 
  The group $\mathrm{Aut}_{B}P$ naturally has a left action on $P$; take its trivial action on
  $E_GO(n)$.
  The group $G$ acts on $P$ and $E_GO(n)$ because they are $G$-vector bundles.
  One can check by \autoref{rem:semidirect-group} that this gives a $\Gamma$-action on $P$ and
  $E_GO(n)$, thus by conjugation on $\mathrm{Hom}(P, E_GO(n))$.   Explicitly,
  the action is
\begin{equation*}
  \begin{array}{ccc}
    (\mathrm{Aut}_{B}P \rtimes G) \times \mathrm{Hom}(P, E_GO(n)) & \to & \mathrm{Hom}(P, E_GO(n)) \\
    ((s, g) , f) & \mapsto & g f s^{-1}g^{-1}.
  \end{array}
\end{equation*}
Since $s \in \mathrm{Aut}_BP$ restricts to identity on $B$, we have
$$ \pi(gfs^{-1}g^{-1}) = g\pi(f)g^{-1}.$$
By \autoref{defn:GGammaprincipal}, the map $\pi$ is a principal $(\mathrm{Aut}_{B}P; \Gamma)$-bundle.

It remains to show that
\begin{equation*}
\mathrm{Hom}(P,E_GO(n))^{\Lambda} \simeq * \text{ for any } \Lambda \subgroup  \Gamma \text{ such
  that } \Lambda \cap \mathrm{Aut}_BP = e.
\end{equation*}
Such subgroups $\Lambda$ are isomorphic to subgroups $H$ of $G$.
The claim follows from various applications of the postponed 
\autoref{lem:contractible}, and it is essentially a consequence of the universality of $E_GO(n)$.

To see it, we first consider the case $\Lambda = H$, that is, the case where $\rho(h) = e$ for all $h \in H$ in
\autoref{rem:Lambeda}.
By restricting the $G$-action to an $H$-action, $E_GO(n)$ is also the universal principal $H$-$O(n)$-bundle.
Then $\mathrm{Hom}(P,E_GO(n))^{H} \simeq *$ by taking $\Pi = O(n)$, $G = H$ and
$\Gamma=O(n) \times H$ in \autoref{lem:contractible}.

In the general case, $\Lambda$ is isomorphic to a subgroup $H \subgroup G$ by the projection map
$\Gamma \to G$, with a possibly non-trivial map $\rho$ in \autoref{rem:Lambeda}.
Here is the crux: the elements in $\mathrm{Aut}_BP$ are $O(n)$-equivariant maps, so the
$(\Gamma = \mathrm{Aut}_BP \rtimes G)$-action on $P$
defined at the beginning of this proof commutes with the
$O(n)$-action; and we have $\Lambda \subgroup \Gamma$.
In other words, $P$ is also a principal $\Lambda$-$O(n)$-bundle.
Since $\Lambda$ acts by $H$ on $E_GO(n)$,  the space $E_GO(n)$ is also the universal principal
$\Lambda$-$O(n)$-bundle.
Now we are basically in the first case again: $\mathrm{Hom}(P,E_GO(n))^{\Lambda} \simeq *$ by taking $\Pi =
O(n)$, $G = \Lambda$ and $\Gamma=O(n) \times \Lambda$ in \autoref{lem:contractible}.
\end{proof}

The following lemma is a consequence of the universality: 
\begin{lem}
  \label{lem:contractible} Let $1 \to \Pi \to \Gamma \to G \to 1$ be an extension of groups.
  Let $$p_{\Pi;\Gamma}: E(\Pi;\Gamma) \to B(\Pi;\Gamma)$$ be the universal principal
  $(\Pi;\Gamma)$-bundle and let
  $p: P \to B$ be any principal $(\Pi;\Gamma)$-bundle. Then $\big(\mathrm{Hom}(P,E(\Pi;\Gamma))\big)^G$ is contractible.
\end{lem}
\begin{proof}
  To clarify the notations, $\mathrm{Hom}(P, E(\Pi;\Gamma))$ is the space of maps of
  (nonequivariant) principal $\Pi$-bundles. By definition,
\begin{equation*}
\mathrm{Hom}(P, E(\Pi;\Gamma)) \cong \mathrm{Map}_{\Pi}(P, E(\Pi;\Gamma)).
\end{equation*}
The space $\mathrm{Hom}(P, E(\Pi;\Gamma))$ has a $\Gamma$-action by
  conjugation. Since $\Pi \subgroup  \Gamma$ acts trivially, it descends
  to a $G$-action, and 
\begin{equation*}
\big(\mathrm{Hom}(P, E(\Pi;\Gamma))\big)^G \cong \mathrm{Map}_{\Gamma}(P, E(\Pi;\Gamma)).
\end{equation*}

By definition, the space $\mathrm{Map}_\Gamma(P, E(\Pi;\Gamma))$ is in fact the space of morphisms
between principal $(\Pi;\Gamma)$-bundles.
  It is non-empty because it consists of the classifying map of $p$.
  It is further path-connected because any two $\Gamma$-maps $P \to E(\Pi;\Gamma)$ will
  both restrict to a classifying map $B \to B(\Pi;\Gamma)$ of $p$, so they are $G$-homotopic.
  The pull back of $p_{\Pi;\Gamma}$ along this homotopy gives a homotopy, or path, between the two
  maps.

  Using the arbitrariness of $P$ in the above argument, one can further show that the space
  $\mathrm{Map}_\Gamma(P, E(\Pi;\Gamma))$ is contractible as follows.
Let $Y$ be a random $G$-space. We denote by $Y \times P$ the principal $(\Pi;\Gamma)$-bundle $Y \times P \to
Y \times B$. Here, $\Gamma$ acts on $Y$ by pulling back the $G$-action and acts $Y \times P$
diagonally. Then we have an adjunction:
\begin{equation}
  \label{eq:2}
\mathrm{Map}_{G}(Y,  \mathrm{Hom}(P, E(\Pi;\Gamma)) )\cong \mathrm{Map}_{\Gamma}(Y \times P,  E(\Pi;\Gamma)).
\end{equation}
By what has been shown, the right hand side, thus the left hand side
of \autoref{eq:2} is always non-empty and path-connected
for any $Y$. Taking
$Y = \mathrm{Hom}(P, E(\Pi;\Gamma))$, we obtain that
$\mathrm{Map}_G(Y,Y)$ is path-connected. In particular, the identity map and the
constant map to a point in $Y^G$ are connected by a path. This implies the contractibility of
$Y^{G} = \big( \mathrm{Hom}(P, E(\Pi;\Gamma))\big)^G$.
\end{proof}

\begin{rem}
  Alternatively, one can show $\mathrm{Map}_{\Gamma}(P,  E(\Pi;\Gamma)) \simeq *$
  using the fact that $E(\Pi;\Gamma)$ is a universal space for a family of
  subgroups of $\Gamma$ specified by \autoref{thm:universalbundle},
  which contains all the isotropy groups of $P$. 
\end{rem}

\subsection{Free loop spaces and adjoint bundles}
\label{sec:free-loop-spaces}
  We end this section by showing an equivariant equivalence of
  the free loop space $LB_G\Pi$ and the adjoint bundle $\mathrm{Ad}(E_G\Pi)$ in \autoref{thm:freeLoop}.
  This gives \autoref{cor:monoidmap}, which upgrades the $G$-equivalence
  $\Omega_bB_GO(n) \simeq O(V)$ to a multiplicative one.
  Our proof follows the non-equivariant treatment in the appendix of Gruher's thesis \cite{Gruher}
  and the key equivariant tool is \autoref{lem:comparefib}.

  We start with $G$-fibrations.
  \begin{defn}
    A $G$-map $p: E \to B$ between $G$-spaces is a $G$-fibration if for all subgroups $H\subgroup G$,
    the map $p^H: E^H \to B^H$ is a Hurewicz fibration.
  \end{defn}

  The first examples of $G$-fibrations are $G$-fiber bundles.
    \begin{exmp}
      Let $p: E \to B$ be a principal $G$-$\Pi$-bundle as in \autoref{defn:Gprincipal}.
      Then $p$ is also a $G$-fibration by
      \autoref{thm:LM}~\autoref{item:LM4}. However, $p: E^H \to B^H$ is not necessarily
      surjective. In contrast to the other parts of \autoref{thm:LM}, we do not have control over
      the components of $B^H$ that are not hit by $p(E^H)$, at least not obviously.
      In this sense, the notion of a $G$-fibration is not as rich as a principal $G$-$\Pi$-bundle.
    \end{exmp}

    \begin{exmp}
      Let $F$ be an effective $\Pi$-space and $q: E' \to B'$ be
      a $G$-fiber bundle with fiber $F$, structure group $\Pi$ as in \autoref{defn:Gvector}. Then
      $q$ is also a $G$-fibration. 
    \end{exmp}

  \begin{lem}\label{lem:Gfibration-fiber}
    We have the following results on the fiber of a $G$-fibration:
\begin{enumerate}
\item Let $p: E \to B$ be a $G$-fibration and $b \in B^H$ be an $H$-fixed point, then the maps
  $(p^{-1}(b))^H \to  E^H \to B^H$ form a fiber sequence.
\item Let $p:D \to B$ and $q: E \to B$ be two $G$-fibrations and $f: D \to E$ be a $G$-map over $B$.
  Take an $H$-fixed point $b \in B^H$.
  If $f$ is a $G$-equivalence, then $p^{-1}(b) \to q^{-1}(b)$ is an $H$-equivalence.
\end{enumerate}
  \end{lem}
  \begin{proof}
    Non-equivariantly ($G = \{e\}$),
    this is the fact that a map over $B$ and homotopy equivalence is a homotopy
    equivalence of fibrations over $B$ (See \cite[7.5-7.6]{Concise}). Equivariantly,
    the first claim is immediate from the definition; the second claim reduces to the
    non-equivariant case for each subgroup $H'\subgroup H$.
  \end{proof}

   We adopt the language of fiberwise monoids in \cite[Definition 4.2.1]{Gruher}.
  \begin{defn}
    A $G$-fibration $p:E \to B$ is a $G$-fiberwise monoid if there is a
    unit section map $\eta: B \to E$ and a multiplication map $m: E
    \times_B E \to E$ over $B$, both $G$-equivariant, that satisfy the unital and associativity conditions. In
    other words, $E$ is a monoid in the category of $G$-fibrations over $B$.
  \end{defn}

  We can relax the strict associativity condition and define $G$-fiberwise
    $A_{\infty}$-monoids as well. Let $\mathscr{A}$ be a reduced $A_{\infty}$-operad in
    $\mathrm{Top}$ ($\mathscr{A}(0)=*$).

  \begin{defn}
    A $G$-fibration $p:E \to B$ is a $G$-fiberwise $A_{\infty}$-monoid if it is an
    algebra over $\mathscr{A}$ in the category of $G$-fibrations over $B$. In concrete words,
    there are $G$-equivariant structure maps over $B$ for each $k \geq 0$
\begin{equation*}
\gamma_k: \mathscr{A}(k) \times_{\Sigma_k} \big( \underbrace{E \times_B E \times_B \cdots \times_B E}_{k \text{
    times}} \big) \to E
\end{equation*}
that satisfy the unital, associativity and $\Sigma$-equivariance conditions of an algebra over an
operad. Here, $\mathscr{A}(k)$ is thought to have the trivial $G$-action; for $k=0$, $\gamma_0: B
\to E$ is just a section of $p$.
\end{defn}

  \begin{defn}
    A morphism of $G$-fiberwise $A_{\infty}$-monoids over $B$
    is a morphism of $A_{\infty}$-monoids in the category of
    $G$-fibrations over $B$.
    An equivalence is a morphism and $G$-equivalence on the total space.
  \end{defn}

  By a $G$-monoid, we mean a monoid in $G$-spaces, and similarly for a $G$-$A_{\infty}$-monoid.  
    Notice that the fiber of a $G$-fiberwise ($A_{\infty}$)-monoid over a point $b \in B$ is not a
    $G$-($A_{\infty}$)-monoid. Instead, it is a
    $G_b$-($A_{\infty}$)-monoid, where $G_b = \{g \in G| gb=b\}$ is the isotropy subgroup of $b$.
    A morphism of fiberwise $G$-($A_{\infty}$)-monoids induces a morphism of $G_b$-($A_{\infty}$)-monoids on
    the fibers over $b$; An equivalence induces a $G_b$-equivalence on the fibers
    by \autoref{lem:Gfibration-fiber}.
  
To clarify this notion, we make the following remarks:
\begin{enumerate}
\item A $G$-fiberwise monoid is a $G$-fiberwise $A_{\infty}$-monoid. In this case,
  the unit section map $\eta$ is $\gamma_0$ and the multiplication map $m$ is $\gamma_2$.
\item The relevant examples of $G$-fiberwise $A_{\infty}$-monoids here are mostly
  $G$-fibrations over $B$ whose fibers are some sort of loops. The structure maps come from
  fiberwise-$A_{\infty}$ structure of loop spaces.
  We will abuse terms to refer to the structure maps as the unit map and the multiplication map.
\item A $G$-fiberwise monoid or a $G$-monoid here is not a ``genuinely
  equivariant algebra'' as it does not have $G$-set indexed multiplications.
\end{enumerate}

  \begin{con}
    For a $G$-space $X$, the free loop space $LX =
    X^{S^1}$ is a $G$-fibration over $X$ by evaluating at a base point of $S^1$.
    It is also a $G$-fiberwise $A_{\infty}$-monoid with the unit map given by the constant loop and the
    multiplication map given by the concatenation of loops.
  \end{con}
  \begin{con}
    For a principal $G$-$\Pi$-bundle $E \to B$, the adjoint bundle of $E$ is
    $Ad(E) = E \times_{\Pi} \Pi_{\mathrm{ad}}$. Here, $\Pi_{\mathrm{ad}}$ is the
    space $\Pi$ with
    adjoint action: for elements $\elm \in \Pi$ and $\ele \in \Pi_{\mathrm{ad}}$,
    $\elm$ acts on $\ele$ by $\elm(\ele) = \elm\ele\elm^{-1}$. 
    As defined, $Ad(E)$ is a $G$-fiber bundle over $B$ with fiber $\Pi$, but no longer a principal
    $G$-$\Pi$-bundle. We claim that $Ad(E)$ has the structure of a $G$-fiberwise monoid over $B$.
    First, $Ad(E)$ is the fiberwise automorphism bundle $\mathcal{I}so_B(E,E)$,
    so naturally a fiberwise monoid over $B$. This is the bundle version of the observation that
    for a right $\Pi$-space $S$ homeomorphic to $\Pi$, there is a homeomorphism
\begin{equation*}
  \begin{array}{ccc}
    \mathrm{Aut}_{\Pi}(S) & \cong & S \times_{\Pi}\Pi_{\mathrm{ad}}\\
    f(s)=s\ele  & \leftrightarrow & [(s,\ele )].
  \end{array}
\end{equation*}
Moreover, $Ad(E) \cong \mathcal{I}so_B(E,E)$ as $G$-spaces, where $G$ acts on $Ad(E)$ by
acting on $E$ and on $\mathcal{I}so_B(E,E)$ by conjugation. This breaks down to commuting the
action of $G$ and $\Pi$ on $E$. Just to clarify the notations,
$$\mathrm{Aut}_{B}(E) = \mathrm{Iso}_B(E,E) \cong \mathrm{Section}( \mathcal{I}so_B(E,E)).$$
  \end{con}

  \begin{thm}
    \label{thm:freeLoop}
    Let $G,\Pi$ be compact Lie groups.
    Then there is a $G$-fiberwise $A_\infty$-monoid
    $(\widetilde{P}E_G\Pi)/\Pi$ over $B_G\Pi$ and equivalences as $G$-fiberwise $A_\infty$-monoids over $B_G\Pi$:
\begin{equation*}
  \begin{tikzcd}
    LB_G\Pi &  (\widetilde{P}E_G\Pi)/\Pi \ar[l,"\xi"',"\simeq"] \ar[r,"\psi","\simeq"'] &  Ad(E_G\Pi)
  \end{tikzcd}
\end{equation*}
\end{thm}

\begin{proof}
  We first construct the space and the map
\begin{equation*}
\widetilde{p}: (\widetilde{P}E_G\Pi)/\Pi \to B_G\Pi.
\end{equation*}
  Recall that $p: E_G\Pi \to  B_G\Pi$ is the universal principal $G$-$\Pi$
  bundle. Denote the space of paths in $E_G\Pi$ that start and end in the same fiber over $B_G\Pi$
  to be
\begin{equation*}
\widetilde{P}E_G\Pi= \{\alpha : I \to E_G\Pi \, | \,  p(\alpha(0)) = p(\alpha(1)) \}.
\end{equation*}
Then $\widetilde{P}E_G\Pi$ inherits an $(\Pi \times G)$-action from $E_G\Pi$. The quotient
$(\widetilde{P}E_G\Pi)/\Pi$ is a $G$-space over $B_G\Pi$ by $\widetilde{p}(\alpha)=p(\alpha(0)).$

The map $\widetilde{p}$ has the structure of a $G$-fiberwise $A_\infty$-monoid.
The unit map $\eta$ 
is given by the constant path in the fiber of $p$.
There is only one constant path in each fiber since we have taken quotient of the $\Pi$-action.
The multiplication map $m$ is given as follows: for two classes of paths
$[\alpha],[\beta] \in (\widetilde{P}E_G\Pi)/\Pi$, we may choose
representatives such that $\alpha(1) = \beta(0) $. Let $m([\alpha],[\beta])=[\alpha.\beta]$ be the
concatenation of the paths: 
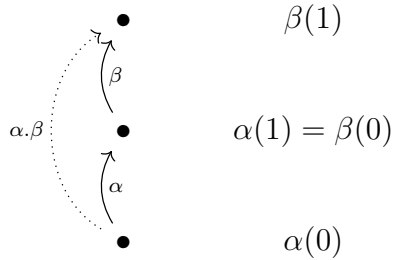
\begin{equation*}
  \begin{tikzcd}
    \bullet & \beta(1)\\
    \bullet \ar[u,bend left, "\beta"']  & \alpha(1) = \beta(0) \\
    \bullet \ar[u,bend left, "\alpha"'] \ar[uu, dotted, bend left = 60, "\alpha.\beta"] & \alpha(0) 
  \end{tikzcd}
\end{equation*}
The class $[\alpha.\beta]$ does not depend on the choice of $\alpha, \beta$.
Both $\eta$ and $m$ are $G$-equivariant.

Next, we compare both $LB_G\Pi$ and $Ad(E_G\Pi)$ with $(\widetilde{P}E_G\Pi)/\Pi$.

On one hand, we have $LB_G\Pi = (\widetilde{P}E_G\Pi)/\Pi^I$. Here, $\Pi^I$ is the group
$\mathrm{Map}([0,1],\Pi)$ and acts on $\widetilde{P}E_G\Pi \subset (E_G\Pi)^I$ pointwise in $I$.
The projection $\widetilde{P}E_G\Pi \to LB_G\Pi$ is a principal $G$-$\Pi^I$-bundle, as the $\Pi^I$
action commutes with the $G$-action on $\widetilde{P}E_G\Pi$.

The projection $\xi: (\widetilde{P}E_G\Pi)/\Pi \to (\widetilde{P}E_G\Pi)/\Pi^I$
commutes with the unit map and multiplication map, so it is a map of $G$-fiberwise $A_\infty$-monoids.
Moreover, we have the following commutative diagram:
\begin{equation*}
  \begin{tikzcd}
    \Pi \ar[r] \ar[d] & \Pi^I \ar[d] \\
    \widetilde{P}E_G\Pi \ar[r,equal] \ar[d] & \widetilde{P}E_G\Pi \ar[d] \\
    (\widetilde{P}E_G\Pi)/\Pi \ar[r,"\xi"] & (\widetilde{P}E_G\Pi)/\Pi^I = LB_G\Pi 
  \end{tikzcd}
\end{equation*}
By \autoref{lem:comparefib}, $\xi$ is a $G$-equivalence. (The idea is that $\Pi$ and $\Pi^I$ are not
so different.)

On the other hand, we may define a $(\Pi \times G)$-equivariant map
\begin{equation*}
  \begin{array}{cccc}
    \bar{\psi}: & \widetilde{P}E_G\Pi &\to & E_G\Pi \times \Pi_{\mathrm{ad}} \\
    &\alpha &\mapsto&  (\alpha(1), \ele)
  \end{array}
\end{equation*}
where $\ele \in \Pi$ is the unique element such that $\alpha(1) = \alpha(0)\ele^{-1}$.
We give $E_G\Pi \times \Pi_{\mathrm{ad}}$ the $G$-action on $E_G\Pi$ and
the diagonal $\Pi$-action.
To check the equivariance of $\bar{\psi}$, take any $(\elm,g) \in \Pi \times G$, then
$(\elm,g) \circ \alpha(t) = g\alpha(t) \elm^{-1}$ for $t \in [0,1]$. So,
\begin{equation*}
\bar{\psi}((\elm,g)\circ \alpha) = (g\alpha(1)\elm^{-1}, \elm\ele\elm^{-1}) =(\elm,g) \circ \bar{\psi}(\alpha). 
\end{equation*}

Since $Ad(E_G\Pi) = (E_G\Pi \times \Pi_{\mathrm{ad}})/\Pi$, we get a map
$\psi: (\widetilde{P}E_G\Pi)/\Pi \to  Ad(E_G\Pi)$. It is easy to check that $\psi$ commutes with the
unit and multiplication maps, and is thus a map of $G$-fiberwise $A_{\infty}$-monoids.

To show that $\psi$ is a $G$-equivalence, we consider the following morphism of principal $G$-$\Pi$-bundles:
\begin{equation*}
  \begin{tikzcd}
    \Pi \ar[r,equal] \ar[d] & \Pi \ar[d] \\
    \widetilde{P}E_G\Pi \ar[r," \bar{\psi}"] \ar[d] & E_G\Pi \times \Pi_{\mathrm{ad}} \ar[d] \\
    (\widetilde{P}E_G\Pi)/\Pi \ar[r,"\psi"] & Ad(E_G\Pi) 
  \end{tikzcd}
\end{equation*}
By \autoref{lem:comparefib}, it suffices to show that $\bar{\psi}$ is a $\Lambda$-equivalence for
any $\Lambda \subgroup  \Pi \times G$ with $\Lambda \cap \Pi = e$.

We can construct a $\Lambda$-homotopy inverse for
$\bar{\psi}: \widetilde{P}E_G\Pi \to {E_G\Pi \times \Pi_{\mathrm{ad}}}$, called $\bar{\phi}$.
The idea is already in Gruher's proof \cite{Gruher}.
But in the equivariant case, $\bar{\phi}$ is dependent on the subgroup $\Lambda$.
(In particular, it is not meant to be a $(\Pi \times G)$-homotopy inverse.)
Recall that $\bar{\psi}$ records the two endpoints of a path.
So an inverse $\bar{\phi}$ is going to choose a canonical path between any two points in a
continuous way.
This choice of canonical path exists because of the $\Lambda$-contractibility of $E_G\Pi$;
it is not meant to be a canonical choice.

The construction of $\bar{\phi}$ is as follows:
Since $E_G\Pi$ is $\Lambda$-contractible, $(E_G\Pi)^{\Lambda}$ is non-empty. We fix a
$\Lambda$-fixed base point $z_0 \in E_G\Pi$. Let $E_G\Pi \times I \to E_G\Pi$ be a
$\Lambda$-equivariant contraction of $E_G\Pi$ to $z_0$; the adjoint of it gives a $\Lambda$-map
$\gamma: E_G\Pi \to P_{z_0}E_G\Pi$. For $z \in E_G\Pi$, we write $\gamma(z)$ as $\gamma_z$.
It is a path connecting $z$ to $z_0$. Now, recall that for an element
$(z,\ele) \in E_G\Pi \times \Pi_{\mathrm{ad}}$, the image $\bar{\phi}(z,\ele) \in
\widetilde{P}E_G\Pi$ wants to be a path from $z \ele$ to $z$  in $E_G\Pi$.
We define it to be
\begin{equation*}
\bar{\phi}(z,\ele) = \text{concatenation of
}\gamma_{z\ele} \text{ and the reverse of } \gamma_{z},
\end{equation*}
 as illustrated in the picture on the left:

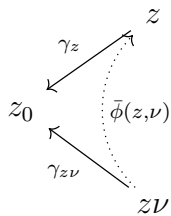
\begin{minipage}{.5\textwidth}
\begin{equation*}
  \begin{tikzcd}
    & z \arrow[dl,"\gamma_{z}"']  \\
    z_0 &  \\
    & z\ele \arrow[ul,"\gamma_{z\ele}"] \arrow[bend left = 45, dotted]{uu}[swap]{\bar{\phi}(z,\ele)}  
  \end{tikzcd}
\end{equation*}
\end{minipage}
\begin{minipage}{.5\textwidth}
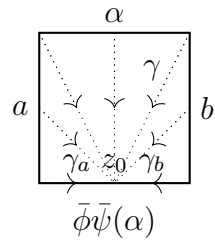
\centering
\begin{tikzpicture}
  \draw [thick]  (0,0) -- (0,2) -- (2,2) -- (2,0) -- cycle;
  \draw [dotted] (0,0)-- (1,0) (2,0) -- (1,0) (0,2) -- (1,0) (2,2) -- (1,0)
  (0,1) -- (1,0) (1,2) -- (1,0) (2,1) -- (1,0);
  \draw [arrows={->[scale = 1.5]} , dotted] (0,0)-- (0.5,0);
  \draw [arrows={->[scale = 1.5]} , dotted] (2,0) -- (1.5,0);
  \draw [arrows={->[scale = 1.5]} , dotted] (0,2) -- (0.5,1);
  \draw [arrows={->[scale = 1.5]} , dotted] (2,2) -- (1.5,1);
  \draw [arrows={->[scale = 1.5]} , dotted] (0,1) -- (0.5,0.5);
  \draw [arrows={->[scale = 1.5]} , dotted] (1,2) -- (1,1);
  \draw [arrows={->[scale = 1.5]} , dotted] (2,1) -- (1.5,0.5);
  \node [above] at (1,2) {$\alpha$};
  \node [above] at (1,0) {$z_{0}$};
  \node [left] at (0,1) {$a$};
  \node [right] at (2,1) {$b$};
  \node [above] at (0.5,0) {\small $\gamma_a$};
  \node [above] at (1.5, 0) {\small $\gamma_b$};
  \node at (1, -0.5) {$\bar{\phi}\bar{\psi}(\alpha)$};
  \node at (1.5,1.5) {$\gamma$};
\end{tikzpicture}
\end{minipage}

It remains to verify that $\bar{\phi}$ is $\Lambda$-homotopy inverse of $\bar{\psi}$.
It is clear that $\bar{\psi}\bar{\phi} = \mathrm{id}$. The illustration above on the
right shows how a $\Lambda$-equivariant homotopy $\bar{\phi}\bar{\psi} \simeq \mathrm{id}$ is defined:
For a path $\alpha \in \widetilde{P}E_G\Pi$ going from a point $a$ to a point $b$, the path
$\bar{\phi}\bar{\psi}(\alpha)$ is the concatenation of $\gamma_{a}$ and the reverse of $\gamma_{b}$.
A homotopy of paths $\bar{\phi}\bar{\psi}(\alpha) \simeq \alpha$ is a map $H$ out of the square, such
that the value of $H$ has been given on the border as indicated.
To fill the interior, we connect every point $x$ on the border to the point labeled by $z_0$ with
line segments and use the map $\gamma_{H(x)}$ on each segment. This homotopy $H$ is ``functorial'' for
elements $\alpha \in \widetilde{P}E_G\Pi$, so it extends to a homotopy $\bar{\phi}\bar{\psi} \simeq
\mathrm{id}$; it is $\Lambda$-equivariant because the map $\gamma$ is.
\end{proof}

We review the Moore loop space construction. For any space $X$ and base point
$b$, the Moore loop space of $X$ at the base point $b$,  $\moore_bX$, is defined to be 
\begin{equation*}
\moore_bX = \{(l,\alpha) \in \bR_{\geq 0} \times X^{\bR_{\geq 0}} | \alpha(0) = b, \ \
\alpha(t)=b \text{ for }t \geq l\}.
\end{equation*}
It has the same homotopy type as $\Omega_bX$ and it is a (strictly associative) monoid, with $\eta: * \to \moore_bX$ given by $\eta(*)= (0,b)$ and $m : \moore_bX \times \moore_bX \to \moore_bX$ given by
$$m((l, \alpha), (s, \beta)) = (l+s, \alpha.\beta) \text{ for } (\alpha.\beta)(t) = \begin{cases}
  \alpha(t) & t < l \\
  \beta(t-l) & l \leq t < l+s \\
  b & t \geq l+s
\end{cases}.$$

As a corollary of \autoref{thm:freeLoop}, we can upgrade \autoref{thm:BGO(n)2}~\autoref{item:BGO(n)loop} into an equivalence
of $G$-$A_{\infty}$-monoids $\Omega_bB_GO(n) \simeq O(V)$.  Strictifying $\Omega_bB_GO(n)$ to the Moore
loop space $\moore_b B_GO(n)$, there is an equivalence of $G$-monoids $\moore_b
B_GO(n) \simeq O(V)$:
\begin{cor}
  \label{cor:monoidmap}
  Take a $G$-fixed base point $b \in B_GO(n)$ in the component indexed by $V$.
  Then $\moore_b B_GO(n)$ is equivalent to $O(V)$ as a $G$-monoid.
  Here, $G$ acts on $\moore_b B_GO(n)$ by acting on $B_GO(n)$ and acts on $O(V)$
  by conjugation.
\end{cor}
\begin{proof}
  We explain how the $G$-$A_{\infty}$-monoid statement is a corollary.
  Take the fiber over $b$ in \autoref{thm:freeLoop} for $\Pi = O(n)$.
  Then there are equivalences of the fibers as $G$-$A_{\infty}$-monoids by \autoref{lem:Gfibration-fiber}.
  The fiber of $LB_GO(n)$ is $\Omega_bB_GO(n)$.
  By  \autoref{thm:BGO(n)2}~\autoref{item:BGO(n)fiber}, the
  fiber of $Ad(E_GO(n))$ is $O(\mathbb{R}^n, V) \times_{O(n)}
  O(n)_{\mathrm{ad}} \cong O(V)$ as $G$-monoid. So there is a zig-zag of
  equivalences of  $G$-$A_{\infty}$-monoids between $\Omega_b B_GO(n) $ and $O(V)$.
  For the $G$-monoid statement, just replace the free loop space and path space in \autoref{thm:freeLoop}
  by the Moore version, and the proof stays intact.

  Explicitly, the zigzag of $G$-monoids is given by 
  \begin{equation}
    \label{eq:OV}
      \begin{tikzcd}
    \moore_b B_GO(n) &  (\widetilde{\moore}_bE_GO(n))/\Pi \ar[l,"\xi"',"\simeq"] \ar[r,"\psi","\simeq"'] &  
 O(V).
  \end{tikzcd}
\end{equation}
  We use $p$ to denote the universal principal $G$-$O(n)$-bundle $E_GO(n) \to B_GO(n)$. We define
\begin{equation*}
\widetilde{\moore}_bE_GO(n)= \{(l,\alpha) | l \in \bR_{\geq 0}, \alpha  : \bR_{\geq 0} \to E_GO(n),
 p(\alpha(0)) = p(\alpha(t)) = b \text{ for }t \geq l\},
\end{equation*}
so that $(\widetilde{\moore}_bE_GO(n))/\Pi = [l, \alpha]$ where the equivalence relation is
\begin{equation*}
  (l,\alpha) \sim (l,\beta) \text{ if there is } \ele \in O(n) \text{ such that }
  \alpha(t) = \beta(t)\ele \text{ for all }t \geq 0.
\end{equation*}
While $\widetilde{\moore}_bE_GO(n) $ does not have the structure of a $G$-monoid,
$(\widetilde{\moore}_bE_GO(n))/\Pi$ does.

Fix a base point $z \in p^{-1}(b) \subset E_GO(n)$.
The maps are given by
\begin{align*}
  \xi([l,\alpha]) & = (l, p(\alpha)) \in \moore_b B_GO(n);\\
  \psi([l,\alpha]) & \in O(V) \text{ is determined by } \alpha(0)\psi([l,\alpha])=\alpha(l). \qedhere
\end{align*}
\end{proof}

\bibliographystyle{alpha}
\bibliography{factorization}

\end{document}